\documentclass[12pt]{article}
\usepackage{mathrsfs}
\usepackage{amssymb}
\usepackage{latexsym}
\usepackage{amsfonts}
\usepackage{amsmath}
\usepackage{eucal}
\usepackage{bm}
\usepackage{bbm}
\usepackage{graphicx}
\usepackage[english]{varioref}
\usepackage[nice]{nicefrac}
\usepackage[all]{xy}
\usepackage{amsthm}
\usepackage{amssymb,amsthm,upref,amscd}

\begin{document}
\renewcommand{\baselinestretch}{1.2}
\setlength{\parindent}{20pt}
\baselineskip18pt
\def\fk{\mathfrak K}

\centerline {\large Irreducibility of Certain Subquotients of}
\centerline {\large Spherical Principal Series Representations of}
\centerline{\large Reductive Groups with Frobenius Maps}

\bigskip

\centerline{Junbin Dong}

\begin{abstract}
For infinite reductive groups with Frobenius maps, we consider the abstract infinite dimensional representations, especially the representations of the groups induced from 1-dimensional representations of Borel subgroups or general parabolic groups. We show that certain subquotients of these induced modules are irreducible. This paper also gives another description of the infinite dimensional Steinberg module. Based on these results, we give a  conjecture about the composition factors of the induced modules.
\end{abstract}


\ \ \ \ \ \ \ \ \ \ \ \ \ \ \ induced modules, Steinberg module.

\section{Introduction}

N. Xi studied abstract representations of infinite reductive groups
with Frobenius maps in [X]. More precisely, he studied modules of
group algebras of these groups and constructed induced modules in
the way for finite groups. It turns out abstract representations of
infinite reductive groups are interesting and many of them are
closely related the representations of finite reductive groups. In
particular, he showed that the infinite dimensional Steinberg module
is irreducible if the ground field of the module is of
characteristic 0 or of the characteristic of the defining field of
the concerned reductive group, by using the irreducibility of the
corresponding Steinberg modules of finite reductive groups. Later,
R. Yang showed that the infinite dimensional Steinberg module is
irreducible for other fields, so that the Steinberg module is always irreducible (see [Y, Theorem 2.2]).

In this paper we are concerned with the composition factors of
representations induced from trivial representation of a Borel
subgroup. We show that certain subquotients constructed in [X, 2.6]
are irreducible when $G$ is of type A or of rank 2, see Theorem 4.1 and Theorem 3.1.
We also give a description of the induced modules of general parabolic groups, see Theorem 6.3.
Analogous to the case of finite reductive group, the infinite dimensional Steinberg module of $\fk G$ can also be expressed as the alternating sum of induced modules of parabolic subgroups in the Grothendieck group, see Theorem 7.1. Based on this, we give a conjecture about the composition factors of
representations induced from trivial representation of a Borel
subgroup, see Conjecture 8.1.

\medskip
The paper is organized as follows. Section 2 contains some
preliminaries, Section 3 deals with rank 2 cases, Section 4 and
Section 5 deal with type A, Section 6 and Section 7 give another description of induced modules and Steinberg module. In section 8  we formulate a
conjecture about the composition factors of representations induced
from trivial representation of a Borel subgroup.

\bigskip

\section{Preliminaries}

In this section we collect some known facts and also establish a few auxiliary results.

{\bf 2.1\ \ } First we recall some basic facts on reductive group defined over a finite field,
one is referred [C] for more details.

Let $G$ be a connected reductive group over the algebraic closure $\bar{\mathbb F}_q$ of a
finite field  $\mathbb{F}_q$ of $q$ elements. We assume that $G$ is defined over $\mathbb{F}_q$.
Then $G$ has a Borel subgroup $B$ defined over $\mathbb{F}_q$ and $B$ contains a maximal torus $T$
defined over $\mathbb{F}_q$. The unipotent radical $U$ of $B$ is also defined over $\mathbb{F}_q$.

For any subgroup $H$ of $G$ defined over $\mathbb{F}_q$ and
any power of $q^a$ of $q$, denote by $H_{q^a}$ the set of
$\mathbb{F}_{q^a}$-points of $H$. Then we have
\begin{equation}
G=\bigcup_{a=1}^{\infty}G_{q^a},\quad B=\bigcup_{a=1}^{\infty}B_{q^a},\quad U=\bigcup_{a=1}^{\infty}U_{q^a},\quad T=\bigcup_{a=1}^{\infty}T_{q^a}.\end{equation}

Let $N=N_G(T)$ be the normalizer of $T$ in $G$. Then $B$ and $N$
forms a $BN$ pairs of $G$. Let $\Phi \subset \text{Hom}(T,\bar{\mathbb F}_q^*)$ be the root system and $\Phi ^+$ be
the set of positive roots determined by $B$.  Let $\Delta$ be the set of simple roots in $\Phi$. Let $W=N/T$ be the Weyl group of $G$ and $S$ the set of simple reflections of $W$.

For each $\alpha\in \Phi$, there is a unique unipotent subgroup
$U_\alpha$ of $G$ which is isomorphic to $\bar{\mathbb F}_q$ and is
stable under conjugation by all elements in $T$. We may choose the
isomorphism $\varepsilon_\alpha:\bar{\mathbb F}_q\to U_\alpha$ so
that $t\varepsilon_\alpha(c)t^{-1}=\varepsilon_\alpha(\alpha(t)c)$.
When $\alpha$ is positive, $U_\alpha$ is in $U$. The unipotent
subgroups $U_\alpha$ are defined over $\mathbb F_q$ and their
$\mathbb F_{q^{a}}$-points are denoted by $U_{\alpha,q^a}$
respectively. The following property is well known.

\medskip

(a) For $w\in W$ and $\alpha\in \Phi$ we have $n_wU_\alpha n_w^{-1}=U_{w(\alpha)}$, here $n_w$ is
a representative of $w$ in $N$. If $n_w$ is in $G_{q^a}$, then $n_wU_{\alpha,q^a} n_w^{-1}=U_{w(\alpha),q^a}$.

\medskip

For simple root $\alpha$ we shall denote by $s_\alpha$ the
corresponding simple reflection. Let $w=s_{\alpha_k} \dots
s_{\alpha_2} s_{\alpha_1}$ be a reduced expression of $w$. Set
$\beta_j=s_{\alpha_1}s_{\alpha_2}\dots s_{\alpha_{j-1}}(\alpha_j)$
for $j=1,\dots,k$. Define
\begin{equation} U_w=U_{\beta_k}\dots U_{\beta_2}U_{\beta_1}\qquad U'_w
=\underset{\beta \in \Phi^+,w(\beta)\in \Phi^+}{\prod}
U_{\beta}.\end{equation} Then

\medskip

(b) $U_w$ and $U'_w$ are subgroups and $wU'_ww^{-1} \in U$;

\medskip

(c) $U=U'_wU_w=U_wU'_w$;

\medskip

(d) each $u \in U_w$ is uniquely expressible in the form
$u=u_{\beta_k}\dots u_{\beta_2}u_{\beta_1}$ with $u_{\beta_i} \in
U_{\beta_i}$.

\medskip

Note that for simple reflection $w=s_\alpha$ we have $U_w=U_\alpha$ and
$U'_w=\underset{\beta \in \Phi^+ -\{\alpha\}}{\prod}
U_{\beta}$ which will be denoted by $U'_\alpha$.

\medskip

(e) Given  two positive roots $\alpha$ and $\beta$,  there exist
integers $c^{mn}_{\alpha \beta}$ such that
\begin{equation}[\varepsilon_\alpha(a),\varepsilon_\beta(b)]:=\varepsilon_\alpha(a)\varepsilon_\beta(b)\varepsilon_\alpha(a)^{-1}\varepsilon_\beta(b)^{-1}=
\underset{m,n>0}{\prod} \varepsilon_{m\alpha+n\beta}(c^{mn}_{\alpha \beta}a^mb^n)\end{equation}
for all $a,b\in \bar{\mathbb F}_q$, where the product is over all
integers $m,n>0$ such that $m\alpha+n\beta \in \Phi^{+}$, taken
according to the chosen ordering.

\medskip
(f)  Numbering all positive roots in any order $\gamma_1,\gamma_2,
\dots ,\gamma_r$, then $U=U_{\gamma_1}U_{\gamma_2}\cdots
U_{\gamma_r}$.

\medskip

{\bf Definition 2.2}\ \ Keep the notations in (f). Let $D$ be a subset of $\Delta$. An element $u=u_{ 1}
u_{ 2}\cdots u_{ r}\in
U$,   $u_{i} \in U_{\gamma_i}$ for $i=1,2,...,r$, is
called $D$-regular provided that $u_{i} \ne 1$ (the neutral element of $G$) if and only if  $\gamma_i\in D$. If $D$ contains only one simple root $\alpha$,
 $D$-regular elements are simply called as $\alpha$-regular element.

\medskip

According to the formula (3), the definition of $D$-regular is independent of the numbering of positive roots.
When $D=\Delta$ is the set of all simple roots,
$D$-regular elements are just  regular unipotent elements.
For $w\in W$, let $\Delta_w=\{\alpha\in \Delta\,|\,w(\alpha)\in \Phi^-\}.$

\medskip

{\bf Lemma 2.3}\ \ Let $w\in W$. For any element $u\in U_{w,q^m}$
 $n>m$, there exists an element $y\in U_{w,q^n}$ such that $yu\in
U_{w,q^n}$ is $\Delta_w$-regular.

\medskip
Proof\ \ Let $w=s_{\alpha_k} \dots
s_{\alpha_2} s_{\alpha_1}$ be a reduced expression of $w$ and
$\beta_j=s_{\alpha_1}s_{\alpha_2}\dots s_{\alpha_{j-1}}(\alpha_j)$
for $j=1,\dots,k$. Assume that $y=y_{\beta_k}\dots
y_{\beta_2}y_{\beta_1}\in
U_{w,q^n}$, $y_{\beta_i}\in U_{\beta_i}$ for all $i$. If  $y_{\beta_i} \in U_{\beta_i,q^{n}}\backslash
U_{\beta_i,q^m}$ for all $i$, then $y$ is $\Delta_w$-regular. Using   the formula (3) in subsection 2.1
 we see that $yu$ is $\Delta_w$-regular.\qed

\medskip

{\bf 2.4}\ \ Let $\mathfrak K$ be a field. For a group $H$ denote by
$\mathfrak K H$ the group algebra of $H$ over $\mathfrak K$. For a
one-dimensional representation $\theta$ of $T$ over $\mathfrak K$,
let $\mathfrak K_\theta$ be the corresponding $\mathfrak KT$-module,
which will be regarded as $\mathfrak KB$-module through the natural
homomorphism $B\to T$. We define the $\mathfrak KG$ module
$$M(\theta)=\mathfrak K G\otimes_{\mathfrak K B}\mathfrak K_\theta,$$
which is called a {\sl spherical principal series representation of
$G$}. When $\theta$ is the trivial representation of $T$ over
$\mathfrak K$, we write $M(tr)$ for $M(\theta)$ and choose a nonzero
element $1_{tr}$ in $\mathfrak K_{tr}$. For $x$ in $\mathfrak KG$, we simply
denote $x\otimes 1_{tr}$ in $M(tr)$ by $x1_{tr}$. For any element
$t\in T$ and $n\in N$, we have $nt1_{tr}=n1_{tr}$. So
$w1_{tr}=n_w1_{tr}$ is well defined for any $w\in W$. Using the
Bruhat decomposition of $G$ we get the following result.
\begin{equation} M(tr)=\sum_{w\in W}\mathfrak KUw1_{tr}.\end{equation}

For any subset $J$ of $S$, we shall denote $W_J$ be the subgroup of
$W$ generated by $J$ and let $w_J$ be the longest element of $W_J$.
Set $$\eta_J =\underset{w\in W_J}{\sum}(-1)^{l(w)}w1_{tr},$$ where
$l(w)$ is the length of $w$. The following result is proved in [X, Prop. 2.3].

\medskip

(a) $M(tr)_J:=\mathfrak KUW\eta_J$ is a submodule of $M(tr)$. In particular
$St=\mathfrak KU\eta_S$ is a submodule of $M(tr)$ and is called the Steinberg
module of $G$.

\medskip

For any integer $a$, the $\fk G_{q^a}$-module $St_a=\fk U_{q^a}\eta$ is isomorphic to the  ordinary Steinberg module which is first constructed by R.Steinberg in [S].

Let $w=s_\alpha$ be a simple reflection corresponding to the simple
root $\alpha$. Recall that $U_\alpha$ and $U'_\alpha$ stand for
$U_w$ and $U'_w$ respectively. Set $n_\alpha=n_{s_\alpha}$. Assume
that $u \in U_\alpha\backslash\{1\}$ and $h \in W$ such that
$l(hw_J)=l(h)+l(w_J)$. The following result are well known or
established  in the proof of [X, Prop. 2.3].

\medskip

(b) There exists $x,y \in U_\alpha \backslash \{1\} $ and $ t \in T$
such that $n_\alpha u {n_\alpha}^{-1} = x n_\alpha t y .$

\smallskip

(c) If $hw_J\le s_\alpha hw_J$, then we have $n_\alpha u
h\eta_J=n_\alpha h\eta_J\in \mathfrak K UW\eta_J$.

\smallskip

(d) If $s_\alpha h\le h$, then $n_\alpha uh\eta_J = xh\eta_J$, where
$x$ is defined in  (b).

\smallskip

(e) If $h\le s_\alpha h$ but $s_\alpha hw_J\le hw_J$, then $n_\alpha
uh\eta_J = (x-1)h \eta_J$, where $x$ is defined in  (b).

\medskip

{\bf 2.5} Let $J$ be a subset of $S$ and $M(tr)'_J$ be the sum of
all $M(tr)_K$ with $J\subsetneq K$. Then $M(tr)'_J$ is a proper
submodule of $M(tr)_J$. Following [X, 2.6], we define $$E_J= M(tr)_J/M(tr)'_J.$$ According to
[X, Prop. 2.7] we have the following result.

\medskip

(a) If $J$ and $K$ are different subsets of $S$ then $E_J$ and $E_K$ are
not isomorphic.

\medskip

For $w\in W$, set $R(w)=\{s\in S\mid ws< w\}$. For any subset
$J$ of $S$, define
\begin{align*} X_J&=\{x\in W\mid x \ \text{has minimal
length in}\  xW_J\},\\
Y_J&=\{w\in X_J\mid R(ww_J)=
J\}, \end{align*}
where $w_J$ is the longest element in the parabolic subgroup
$W_J$.

\def\fk{\mathfrak K}

\def\st{\stackrel}
\def\sc{\scriptstyle}

For any $w\in W$, let
$$C_w=\sum_{y\le w}(-1)^{l(w)-l(y)}p_{y,w}(1)y\in\mathfrak KW,$$
 where $p_{y,w}$ are  Kazhdan-Lusztig polynomials. Then the elements $C_w,\ w\in W$, form a basis of $\mathfrak KW,$ see [KL].

 \medskip

 {\bf Lemma 2.6} For any subsets $J,\ K$ of $S$, we have

 (a) the elements $xC_{w_K}, \ x\in X_K$, form a basis of $\fk WC_{w_K}$;

 (b) the elements $C_{xw_K},\ x\in X_K$, form a basis of $\fk WC_{w_K}$;

 (c) the elements $wC_{w_J}, \ C_{xw_J}, \ w\in Y_J,\ x\in X_J\backslash Y_J$, form a basis of  $\fk WC_{w_J}$.

 Proof. Since $W=X_KW_K$ and $wC_{w_K}=(-1)^{l(w)}C_{w_K}$ for any $w\in W_K$, we see that (a) is true.
By Lemma 2.8 (c) in [G], for $x\in X_K$, we have
\begin{equation}C_{xw_K}=xC_{w_K}+ \sum_{\st{\sc y\in X_K}{y<x}}a_yyC_{w_K},\quad a_y\in\fk.\end{equation}
Using induction on $l(x)$ we see that
\begin{equation}xC_{w_K}=C_{xw_K}+ \sum_{\st{\sc y\in X_K}{y<x}}a'_yC_{yw_K},\quad a'_y\in\fk.\end{equation}
(b) follows.

We claim that for any $x\in X_J$, the element $xC_{w_J}$ is a linear combination of the elements $wC_{w_J}, \ C_{xw_J}, \ w\in Y_J,\ x\in X_J\backslash Y_J$.
If $l(x)=0,$ then $x\in Y_J$, the claim is true. Now assume that the claim is true for $y\in X_J$ with $l(y)<l(x)$. If $x$ is in $Y_J$, the claim is clear.
If $x$ is in $X_J\backslash Y_J$, using formula (6) and induction hypothesis we see that the claim is true. (c) is proved.\qed

\medskip
{\bf Lemma 2.7} For any subset $J\subseteq S$, denote by $C_J$ the
image of $\eta_J$ in $E_J$. Then the $\mathfrak KG$-module $E_J$ is
the sum of all $\mathfrak KU_{w_Jw^{-1}}\,w C_J$, $w\in Y_J$, i.e.,
$$E_J=\sum_{w\in Y_J}\mathfrak KUw C_J=\sum_{w\in Y_J}\mathfrak KU_{w_Jw^{-1}}\,w C_J.$$

Proof. For any $w\in W$,  set $$c_w=(-1)^{l(w)}C_w1_{tr}\in M(tr).$$
  By
Lemma 2.6 (vi) in [KL], we  have $$c_{w_K} = \eta_K\quad \text{for any subset $K$ of } S.$$
According Lemma 2.6 (c), we have
\begin{equation}M(tr)_J=\fk UW\eta_J=\sum_{w\in Y_J}\fk Uw\eta_J+\sum_{x\in X_J\backslash Y_J}\fk Uc_{xw_J}.\end{equation}

For $x\in  X_J\backslash Y_J$, we have $J\subsetneq R(xw_J)=K$. Then
$xw_J=x'w_K$ for some $x'\in X_K$. Note that $c_{w_K}=\eta_K$. By
Lemma 2.6 (b), $c_{xw_J}=c_{x'w_K}=h\eta_K$ for some $h\in \fk W$.
Therefore, $c_{xw_J}$ is in $M(tr)'_J$ if $x\in  X_J\backslash Y_J$.
Using formula (7) we get
$$E_J=M(tr)_J/M(tr)'_J=\sum_{w\in Y_J}\mathfrak KUw C_J.$$

To finish the proof it suffices to prove that
$Uw\eta_J=U_{w_Jw^{-1}}\,w\eta_J$ for $w\in Y_J$. Since
$\eta_J=\sum_{z\in W_J}(-1)^{l(z)}z1_{tr}$, we only need to show
that $Uwz1_{tr}=U_{w_Jw^{-1}}\,wz1_{tr}$ for any $z\in W_J$ and
$w\in Y_J$.

Assume that $\alpha$ is a positive root in $\Phi$. For any $y\in W_J$ we prove that $yw^{-1}(\alpha)\in \Phi^+$ if $w_Jw^{-1}(\alpha)$ is positive.
Otherwise,  $yw^{-1}(\alpha)=\beta$ was negative root. Since $w_J=w_Jy^{-1}y$ and $l(w_Jy^{-1})+l(y)=l(w_J)$, we have $\gamma=w_Jy^{-1}(\beta)=w_Jw^{-1}(\alpha)\in \Phi^+$.
This forces that $\beta$ is a negative root in the root system $\Phi_J$ corresponding to $W_J$ and $\gamma\in \Phi_J\cap \Phi^+$. Then $\delta=w_J(\gamma)$ is a negative root
in $\Phi_J$. Since $l(ww_J)=l(w)+l(w_J)$, we have $\alpha=ww_J(\gamma)=w(\delta)\in \Phi^-$. This contradicts the assumption $\alpha$ being positive.
Therefore we have  $z^{-1}w^{-1} U'_{w_Jw^{-1}}wz\in U$. Since $U1_{tr}=1_{tr}$, we see $U'_{w_Jw^{-1}}wz1_{tr}=wz1_{tr}$ for $w\in Y_J$ and $z\in W_J$.
Using 2.1(c) we get that $Uwz1_{tr}=U_{w_Jw^{-1}}\,wz1_{tr}$ for any $z\in W_J$ and $w\in Y_J$. The lemma is proved.\qed

\medskip

{\bf 2.8} An interesting question is whether these $\mathfrak
KG$-module $E_J$ are irreducible. When $J$ is the empty set, $E_J$
is just the trivial representation of $G$. When $J=S$, $E_J$ is just
the Steinberg module, which is irreducible (see [X, Theorem 3.2]
and [Y, Theorem 2.2]).

Provided  that char \,$\mathfrak K \ne$\,char\,$\mathbb F_q$, we
will show that $E_J$ is irreducible in the following cases: (1)
$G$ is of rank 2 (see Theorem 3.1), (2)  $G$ is  of type A and $J$
contains only one element or $J$ is a maximal proper subset of $S$
(see Theorem 4.1).

The following result will be used frequently in our proof for Theorems 3.1 and 4.1.

\medskip

{\bf Lemma 2.9 } Let $  M$ be a  $\mathfrak KG$-module and $\eta\in
M$ is $T$-fixed (i.e., $t\eta=\eta$ for all $t\in T$). Assume char \,$\mathfrak K \ne$\,char\,$\mathbb F_q$. If $  M'$ is a submodule of $  M$ containing $\underset{x\in U_{q^a}}{\sum}x\eta$
for some positive integer $a$, then $\eta \in   M'$.

Proof. The argument in [Y, 2.7] works well here. Let
$s_{\alpha_r}s_{\alpha_{r-1}}\dots s_{\alpha_1}$ be a reduced
expression of the longest element $w_0$ of $W$ . Set
$\beta_i=s_{\alpha_1}s_{\alpha_{2}}\dots s_{\alpha_{i-1}
}{(\alpha_i)}$. Then for any positive integer $b$,
$$X_{i,q^b}=U_{\beta_r,q^b}U_{\beta_{r-1},q^b}\dots U_{\beta_i,q^b}$$
is a subgroup of $U$. Clearly, $X_{i,q^b}$ is a subgroup of
$X_{i,q^{b'}}$ if $\mathbb{F}_{q^b}$ is a subfield of
$\mathbb{F}_{q^{b'}}$. We understand that $X_{r+1,q^b}=\{1\}.$

First we use induction on $i$ to show that there exists positive
integer $b_i$ such that the element $\underset{x\in
X_{i,q^{b_i}}}{\sum}x\eta$ is in $   M'$. When $i=1$, this is true
for $b_1=a$ by assumption. Now we assume that $\underset{x\in
X_{i,q^{b_i}}}{\sum}x\eta$ is in $   M'$. We show that
$\underset{x\in X_{{i+1},q^{b_{i+1}}}}{\sum}x\eta$ is in $   M'$ for
some $b_{i+1}$.

Let $c_1,c_2,\dots,c_{q^{b_i}+1}$ be a complete set of
representatives of all cosets of $\mathbb{F}_{q^{b_i}}^*$ in
$\mathbb{F}_{q^{2b_i}}^*$. Choose $t_1,t_2,\dots,t_{q^{b_i}+1}\in T$
such that $\beta_i(t_j)=c_j$ for $j=1,2,\dots,q^{b_i}+1$. Note that
$t\eta=\eta$ for any $t\in T$. Thus$$\sum_{j=1}^{q^{b_i}+1}t_j
\underset{x\in U_{\beta_i,q^{b_i}}}{\sum} x\eta =
q^{b_i}\eta+\underset{x\in U_{\beta_i,q^{2b_i}}}{\sum} x\eta.$$
Since $X_{i,q^{b_i}}=X_{i+1,q^{b_i}}U_{\beta_i,q^{b_i}}$ and
$\underset {x\in X_{i,q^{b_i}}}{\sum}x\eta$ is in $   M'$, we see
that
$$\zeta= \sum_{j=1}^{q^{b_i}+1}t_j \underset{x\in X_{i,q^{b_i}}}{\sum}x\eta
= \sum_{j=1}^{q^{b_i}+1}t_j \underset{y\in X_{i+1,q^{b_i}}}{\sum}y
\underset{x\in U_{\beta_i,q^{b_i}}}{\sum}x\eta$$
$$= \underset {y\in
 X_{i+1,q^{b_i}}}{\sum}\sum_{j=1}^{q^{b_i}+1}t_jyt_j^{-1}(t_j\underset{x\in U_{\beta_i,q^{b_i}}}{\sum}x\eta) \in   M'.$$
Choose $b_{i+1}$ such that all $\beta_m(t_j)(r\geq m \geq i)$ are
contained in $\mathbb{F}_{q^{b_{i+1}}}$ and
$\mathbb{F}_{q^{b_{i+1}}}$ contains  $\mathbb{F}_{q^{2b_i}}$. Then
$t_j y t_j^{-1}$ is in $X_{i+1,q^{b_{i+1}}}$ for any $y\in
X_{i+1,q^{b_i}}$. Let $Z\in \mathfrak KG$ be the sum of all elements
in $X_{i+1,q^{b_{i+1}}}$. Then we have
\begin{equation}\displaystyle Z\zeta = q^{(r-i)b_i}Z \sum_{j=1}^{q^{b_i}+1}t_j \underset {x\in
U_{\beta_i,q^{b_i}}}{\sum}x\eta= q^{(r-i)b_i}Z(q^{b_i}\eta +
\underset{x\in U_{\beta_i,q^{2b_i}}}{\sum}x\eta)\in
M'.\end{equation}

Because $\underset {x\in X_{i,q^{b_i}}}{\sum}x\eta$ is in $   M'$,
we have $\underset {x\in X_{i,q^{2b_i}}}{\sum}x\eta\in   M'$. Thus

\begin{equation}Z\underset {x\in
X_{i,q^{2b_i}}}{\sum}x\eta=q^{2(r-i)b_i}Z\underset{x\in
U_{\beta_i,q^{2b_i}}}{\sum}x\eta\in   M'.\end{equation}

Since $q\neq 0$ in $\mathfrak K$, combining formula (8) and (9) we
see that $Z\eta \in   M'$, i.e.,$\underset {x\in
X_{i+1,q^{b_{i+1}}}}{\sum}x\eta$ is in $   M'$. Note that
$X_{r,q^{b_r}}=U_{\beta_r,q^{b_r}}$. Now we have $\underset {x\in
U_{\beta_r,q^{b_r}}}{\sum}x\eta \in   M'$ and $\underset {x\in
U_{\beta_r,q^{2b_r}}}{\sum}x\eta \in   M'$. Applying formula (8) to
the case $i=r$ we get that $q^{b_r}\eta + \underset {x\in
U_{\beta_r,q^{2b_r}}}{\sum}x\eta \in
  M'$. Therefore $\eta$ is in $   M'$. The lemma is proved.\qed

\bigskip

\section{Rank 2 cases}

In this section we consider the irreducibility of the modules $E_J$ for rank 2 cases. The main result is the following.

\medskip

{\bf Theorem 3.1} Assume that $G$ is of rank 2  and char\,$\mathfrak
K \ne\text{char\,}\bar{\mathbb{F}}_q$. Then the $\mathfrak
KG$-modules $E_J$ are always irreducible.

\medskip

We prove Theorem 3.1 case by case. Since $S$ has only two
elements, $E_S$ and $E_\emptyset$ are irreducible,   we only need to consider the case $J$ containing one
element. We shall write $E_d$ instead of $E_J$ when $J=\{d\}$. The neutral element in $W$ will denoted by $e$.

\medskip

{\bf 3.2 }\ In this subsection $G$ is assumed of type $A_2$.
Let $\alpha,\ \beta$ be the simple roots of $\Phi$ and denote by $s,r$
the corresponding simple reflections. Then $W=\{e, \ s,\ r,\ sr,\ rs,\ srs\}$ and $Y_{\{s\}}=\{e,\ r\}$. Since $U_s=U_\alpha$ and $U_{sr}=U_{\alpha+\beta}U_\beta$, by Lemma 2.7 we get
\begin{equation} E_{s}=\mathfrak KUC_s+\mathfrak KUrC_s=\mathfrak KU_\alpha C_s+\mathfrak KU_{\alpha+\beta}U_\beta rC_s,\end{equation}
 where $C_s$ is the image of $(1-s )1_{tr}$
in $E_{s }$.

We will show that $E_s$ is generated by any nonzero element in $E_s$, so that $E_s$ is irreducible.
Let $\xi$ be a nonzero element in $E_s$, by formula (10) we then have
\begin{equation}\xi=\underset{u\in U_\alpha}{\sum}a_uuC_s+\underset{u\in U'_\alpha}{\sum}b_uurC_s,\end{equation}
where $U'_\alpha=U_{\alpha+\beta}U_\beta$; $a_u,b_u\in \mathfrak K$ and only
finitely many of them are nonzero.

Clearly there is a positive integer $m$ such that $u\in U_{q^m}$
whenever $a_u\ne 0$ or $b_u\ne 0$. We may require that $a_u=0$ for
some $u\in U_{\alpha,q^m}$, this can be done by choosing
sufficiently large $m$. We may further require that $a_1=0$ (recall
that 1 also stands for the neutral element of $G$). Otherwise,
choose $y\in U_{\alpha,q^m}\setminus \{1\}$ such that
$a_{y^{-1}}=0$, then $y\xi$ satisfies the requirement and we then
consider $y\xi$ instead of $\xi$.

Let $n_s$ (resp. $n_r$) be a representative of $s$ (resp. $r$) in $N=N_G(T)$. Then $n_s^2$ and $n_r^2$ are in $T$. Since $m$  is sufficiently large, we may require that  $n_s$ and $n_r$ are in $G_{q^m}$.
By 2.1(a) we have

\begin{equation}n_s(U_{\alpha+\beta,q^m}U_{\beta,q^m}) n^{-1}_s=U_{\beta,q^m}U_{\alpha+\beta,q^m}=U_{\alpha+\beta,q^m}U_{\beta,q^m}.\end{equation}
Noting that the image of $(1-s-r+sr+rs-srs)1_{tr}$ in $E_s$ is zero, we get
\begin{equation}n_srC_s=srC_s= rC_s-C_s.\end{equation}
Using formulas (12) and (13) we get the following claim.

 (a) For any $u\in U_{\alpha+\beta,q^m}U_{\beta,q^m}$ there exists $u'\in U_{\alpha+\beta,q^m}U_{\beta,q^m}$ such that
\begin{equation} n_surC_s=u'(rC_s-C_s).\end{equation}

By 2.4(d), we obtain

(b) if $a_u\ne 0$, then $n_suC_s=(z-1)C_s$ for some $z\in U_{\alpha,q^m}$.
(Note that by assumption $a_u\ne 0$ implies that $u\ne 1.)$

Set$$M'=\mathfrak KG\xi,\quad \mathbb{X}=\displaystyle \sum_{x\in
U_{q^m}}x,\quad \phi_\alpha=\displaystyle {\sum_{u\in
U_\alpha}}a_u,\quad \phi'_\alpha={\sum_{u\in U'_\alpha}}b_u,\quad
\phi=\phi_\alpha+\phi'_\alpha.$$ Note that $\mathbb{X} u=
\mathbb{X}$ for any $u\in U_{q^m}$. Using (a) and (b)  we get
\begin{align*} \mathbb{X}\xi - \mathbb{X} n_s\xi&=\phi_\alpha \mathbb{X}C_s+
\phi'_\alpha \mathbb{X} rC_s-\phi'_\alpha \mathbb{X}(rC_s-C_s)\\
&=\phi \mathbb{X}C_s\in M',\end{align*}
 If $\phi\neq 0$
 then $\mathbb{X} C_s \in
M'$. By Lemma 2.9 , $C_s\in M'$. Thus $M'=E_{s}$ in this case.

 Now assume that  $\phi= 0$
 but $a_{u_0}\ne 0$ for some $u_0\in U_{\alpha,q^m}$.
 Note that $n_sC_s=-C_s$. Using (a) and (b)  we see
$$\mathbb{X} \xi - \mathbb{X} n_s{u_0}^{-1}\xi=\phi \mathbb{X} C_s+ a_{u_0}\mathbb{X} C_s=a_{u_0}\mathbb{X} C_s\in M'.$$ Using Lemma 2.9 , we get
$C_s\in M'$. So in this case we also have $M'=E_{s}$. We have proved
the following result.

(c) Let $\xi$ be a nonzero element in $E_s$ of the form (11) with  $a_1= 0$. If $\phi\ne 0$ or $a_{u_0}\ne 0$ for some $u_0\in U_{\alpha}$, then $\mathfrak KG\xi=E_s.$

If all $a_u$ are zero, then there exists $v\in
U_{\alpha+\beta,q^m}U_{\beta,q^m}$ such that $b_{v}\ne 0$. Then
$$n_r{v}^{-1}\xi=\underset{u\in U_\alpha}{\sum}a'_uuC_s+\underset{u\in U'_\alpha}{\sum}b'_uurC_s,$$
where $a'_u,b'_u\in \mathfrak K$ and $a'_1\ne 0$.

Choose $z\in U_\alpha\backslash U_{\alpha,q^m}$, let
 $$\xi'=zn_rv^{-1}\xi=\underset{u\in U_\alpha}{\sum}a''_uuC_s+\underset{u\in U'_\alpha}{\sum}b''_uurC_s.$$
 Since $a'_u\ne 0$ implies that $u\in U_{\alpha,q^m}$, we have  $a''_1= 0.$ By (c) we have $\mathfrak KG\xi'=E_s$. So in this case we also have $M'=\mathfrak KG\xi=\mathfrak KGzn_rv^{-1}\xi=\mathfrak KG\xi'=E_s$.

We have proved that $E_s$ is an irreducible $\mathfrak KG$-module. Since $s,r$
are symmetric, the $\mathfrak KG$-module $E_r$ is also irreducible.

\medskip

{\bf 3.3} In this subsection $G$ is assumed of type $B_2$. Let
$\alpha,\beta$ be the simple roots of $\Phi$ and denote by $s,r$ the
corresponding simple reflection. We assume that $\alpha$ is short,
then $\Phi^+$ consists of $\alpha,\ \beta,\ \alpha+\beta,\
2\alpha+\beta.$  We have $W=\{e,\ s,\ r,\ sr,\ rs,\ srs,\ rsr,\
srsr\}$ and $Y_{\{s\}}=\{e,\ r,\ sr\}$.

Since $U_s=U_\alpha,\  U_{sr}=U_{\alpha+\beta}U_\beta,\
U_{srs}=U_{\alpha+\beta}U_{2\alpha+\beta}U_\alpha$, by Lemma 2.7 we
get
\begin{equation}E_s=\mathfrak KU_\alpha C_s+\mathfrak KU_{\alpha+\beta}U_\beta
rC_s+\mathfrak KU_\alpha
U_{\alpha+\beta}U_{2\alpha+\beta}srC_s,\end{equation} where $C_s$ is
the image of $(1-s )1_{tr}$ in $E_{s}$.

We show that $E_s$ is generated by any nonzero element in $E_s$ so
that $E_s$ is irreducible. Let $\xi$ be a nonzero element of $E_s$,
by formula (15) we have
\begin{equation}\xi=\underset{u\in U_\alpha}{\sum}a_uuC_s+
\underset{u\in V}{\sum}b_uurC_s+\underset{u\in U'_\beta}{\sum}c_uusrC_s,\end{equation}
where $V=U_{\alpha+\beta}U_\beta,$ $U'_\beta=U_\alpha
U_{\alpha+\beta}U_{2\alpha+\beta}=
U_{\alpha+\beta}U_{2\alpha+\beta}U_\alpha,$ $a_u, b_u,c_u\in
\mathfrak K$ and only finitely many of them are nonzero.

Clearly there is a positive integer $m$ such that $u\in U_{q^m}$
whenever $a_u\ne 0$, or $b_u\ne 0$, or $c_u\ne 0$. We can choose $m$
so that $n_s$ and $n_r$ are in $G_{q^m}$. Recall the Definition 2.2,
we may further require that $u$ is $\alpha$-regular (resp.
$\beta$-regular; $\alpha$-regular) if $a_u\ne 0$ (resp. $b_u\ne 0$;
$c_u\ne 0$). Otherwise, according to Lemma 2.3, we can choose some
$y \in U$ such that $y\xi$ satisfies the requirement. Then we
consider $y\xi$ instead of $\xi$.

Let $n_s$ (resp. $n_r$) be a representative of $s$ (resp. $r$) in
$N=N_G(T)$. Then $n_s^2$ and $n_r^2$ are in $T$. Since $m$  is
sufficiently large, we may require that  $n_s$ and $n_r$ are in
$G_{q^m}$. Noting that the image of $(1-s-r+
rs+sr-srs-rsr+rsrs)1_{tr}$ in $E_s$ is zero, we have
\begin{equation} rsrC_s =C_s-rC_s+srC_s. \end{equation}

With the assumption of $\xi$, using 2.1(a), 2.4(c), 2.4(d), 2.4(e)
and formula (17) we have

\medskip

(a) If $a_u\ne 0$, then $n_suC_s=(z-1)C_s$ for some $z\in U_{q^m}$.
If $b_u\ne 0$, then $n_surC_s\in U_{q^m} srC_s$. If $c_u\ne 0$, then
$n_susrC_s\in U_{q^m}  srC_s$.

\medskip
(b) If $a_u\ne 0$, then $n_ruC_s\in U_{q^m}rC_s$. If $b_u\ne 0$ then
$n_rurC_s\in U_{q^m}rC_s$. If $c_u\ne 0$ then
$n_rusrC_s=u'(C_s-rC_s+srC_s)$ for some $u'\in U_{q^m}$.

\medskip

For convenience, set $$M'=\mathfrak KG\xi,\quad
\mathbb{X}=\displaystyle \sum_{x\in U_{q^m}}x,\quad
\phi_a=\displaystyle \underset{u\in U_\alpha}{\sum}a_u, \quad
\phi_b=\displaystyle \underset{u\in V}{\sum}b_u, \quad
\phi_c=\displaystyle \underset{u\in U'_\beta}{\sum}c_u.$$

Note that $\mathbb{X}u=\mathbb{X}$ for any $u\in U_{q^m}$. Using (a)
we get
$$\mathbb{X} n_s\xi=(\phi_b+\phi_c)\mathbb{X} srC_s \in M'.$$

If $\phi_b+\phi_c \neq 0$, using Lemma 2.9 we know that $srC_s\in
M'$. Thus $M'=E_s$ in this case.

Let $\xi_1=n_r\xi$ and by (b) we write $\xi_1$ as
\begin{equation}\xi_1=\underset{u\in U}{\sum}f_uu rC_s+ \underset{u\in U}{\sum}g_uu
rsrC_s,\end{equation} where $\phi_f = \phi_a+\phi_b$.

We can assume that all $u\in U$ appear in $\xi_1$ are
$\beta$-regular. Otherwise, we can choose some $y \in U$ such that
$y\xi_1$ satisfies the requirement by Lemma 2.3. Then we consider
$y\xi_1$ instead of $\xi_1$.With the assumption on $\xi_1$, using
2.4(e), it is easy to see the following result.

\medskip

(c) If $g_u\ne 0$, then $n_sursrC_s=(y-z)rsrC_s$ for some $y,z \in
U$.

\medskip
Using (b) and (c) we get
$$\mathbb{X} n_s\xi_1=(\phi_a+\phi_b) \mathbb{X}srC_s  \in M'.$$
If $\phi_a+\phi_b \neq 0$, then $srC_s\in M'$ by Lemma 2.9 which
implies $M'=E_s$.

Now suppose $\phi_a+\phi_b =0$ then using the formula (18) we get
$$\mathbb{X} \xi_1=\phi_c \mathbb{X} rsrC_s \in M'.$$
If $\phi_c \neq 0$, Also by Lemma 2.9 we have $rsrC_s\in M'$ and
$M'=E_s$.

\medskip

(d) By above computation, if $\phi_a\ne 0$ or $\phi_b\ne 0$ or
$\phi_c\ne 0$ then $M'=E_s$.
\medskip

We just need to deal with the case $\phi_a=\phi_b =\phi_c =0$. In
the following we assume that $\phi_a=\phi_b =\phi_c =0$.

Now assume that $a_{u_1}\ne 0$ for some $u_1 \in U_\alpha$. Note
that $n_sU_{q^m}rC_s\in U_{q^m}srC_s$ and $n_sU_{q^m}srC_s\in
U_{q^m}rC_s$ or $U_{q^m}srC_s$. Then we have
\begin{equation}n_s{u_1}^{-1}\xi=
\underset{u\in U_\alpha}{\sum}a'_uuC_s+\underset{u\in
V}{\sum}b'_uurC_s+\underset{u\in U'_\beta}{\sum}c'_uusrC_s,
\end{equation} where $\phi_{a'} =-a_{u_1}$. By the result in (d), $\mathfrak
KGn_s{u_1}^{-1}\xi=E_s$ which implies $M' =E_s$ in this case.

We compute the form of $n_sn_rurC_s$ for any  $u\in
U_{\alpha+\beta}U_\beta$ and  $n_sn_rusrC_s$ for any  $u \in
U_{\alpha+\beta}U_{2\alpha+\beta}U_\alpha$. The following results
can be obtained from the 2.4(c), 2.4(d), 2.4 (e).

\medskip
(e) If $u=1$ then $n_sn_rurC_s=-C_s$. If $u_{\beta}\ne 1$ then
$n_sn_rurC_s\in UsrC_s$. If $u_{\alpha+\beta}\ne 1$ but
$u_{\beta}=1$ then $n_sn_rurC_s=(x-y)C_s$ for some $x,y\in U$.

\medskip

(f) If $u_{\alpha+\beta}=1$ then $n_sn_rusrC_s=-xC_s+xrC_s-xsrC_s$
for some $x\in U$. If $u_{\alpha+\beta}\ne 1$ then
$n_sn_rusrC_s=(y-z)(C_s+srC_s)$ for some $y,z\in U$.

\medskip

Assume that all $a_u$ are zero but $b_{u_2}\ne 0$ for some $u_2 \in
V$, we write $n_sn_r{u_2}^{-1}\xi$ as
\begin{equation}n_sn_r{u_2}^{-1}\xi= \underset{u\in U_\alpha}{\sum}a''_uuC_s+
\underset{u\in V}{\sum}b''_uurC_s+\underset{u\in
U'_\beta}{\sum}c''_uusrC_s\end{equation} and by (e) and (f) we can
see $\phi_{a''}+\phi_{b''}= -b_{u_2}$ which is nonzero. We can deal
with this case by the argument in (d) and obtain $M'=E_s$.

For the last case that $a_u,b_u$ are all zero but $c_{u_3}\ne 0$ for
some $u_3 \in U'_\beta$, we can easy to see that
$$n_s{u_3}^{-1}\xi= \underset{u\in V}{\sum}b^*_uurC_s+\underset{u\in U'_\beta}{\sum}c^*_uusrC_s$$
where $b^*_1= c_{u_3}$ is nonzero. So in this case we also have
$M'=E_s$.

Now we complete the proof that if $\xi $ is a nonzero element of
$E_s$ then $\mathfrak KG \xi=E_s$. So $E_s$ is irreducible.
Similarly, we can prove $E_r$ is irreducible.

\medskip

{\bf 3.4} In this subsection, we give some a remark on the proof of
type $B_2$ case. Instead of computing the results of (e) and (f) in
subsection 3.4, we give another proof which is also useful later.

Let $\xi$ be a element of the form in (16) satisfy our requirement
as before. We have proved if $\phi_a\ne 0$ or $\phi_b\ne 0$ or
$\phi_c\ne 0$ then $\mathfrak K G \xi=E_s$.

When $a_{u_1}\ne 0$ for some $u_1$, we have $\mathfrak K G \xi=E_s$
by the same reason as subsection 3.3. Now suppose all $a_u$ are zero
but $b_{u_2}\ne 0$ for some $u_2 \in V$. Using 2.4(c), 2.4(d),
2.4(e), we have

\medskip

(a) If $u_\beta =1$, then $n_rurC_s\in UC_s$.

\medskip

(b) If $u_\beta \ne 1$, then $n_rurC_s\in UrC_s$.

\medskip

(c) $n_rusrC_s = xrsrC_s= x(C_s-rC_s+srC_s)$ for some $x\in U$.

\medskip

We consider the element $n_r{u_2}^{-1}\xi$ and write it as
\begin{equation}n_r{u_2}^{-1}\xi=\underset{u\in U_\alpha}{\sum}a'_uuC_s+\underset{u\in
V}{\sum}b'_uurC_s+\underset{u\in U'_\beta}{\sum}c'_uusrC_s.
\end{equation}

If some $a'_u$ is nonzero, then we have $E_s=\mathfrak K
Gn_r{u_2}^{-1}\xi\subseteq \mathfrak K G\xi$ which implies
$\mathfrak KG \xi=E_s$. Otherwise, We may choose the isomorphism
$\varepsilon_\beta:\bar{\mathbb F}_q\to U_\beta$ and we consider the
element $n_r\varepsilon_\beta(c) \xi$ and write it as
\begin{equation}n_r\varepsilon_\beta(c) \xi =\underset{u\in
U_\alpha}{\sum}a'_u(c)u(c)C_s+\underset{u\in
V}{\sum}b'_u(c)u(c)rC_s+\underset{u\in
U'_\beta}{\sum}c'_u(c)u(c)srC_s. \end{equation} The number of $c\in
\bar{\mathbb F}_q$ which makes $\varepsilon_\beta(c)u_\beta =1$ for
$b_u\ne 0$ is finite. Since $\underset{u\in
U_\alpha}{\sum}a'_u(c)u(c)C_s$ are not always zero and the roots of
$\underset{u\in U_\alpha}{\sum}a'_u(c)u(c)C_s=0$ is finite, we can
choose some $c_1\in \bar{\mathbb F}_q$ such that $a'_u(c_1)\ne 0$
for some $u$. Therefore we can prove $\mathfrak K G \xi=E_s$.

For the last case that $a_u,b_u$ are all zero but $c_{u_3}\ne 0$ for
some $u_3 \in U'_\beta$, the proof is similar to subsection 3.3 and
we can prove $E_s$ is irreducible.

\medskip

{\bf 3.5} In this subsection $G$ is assumed of type $G_2$. Let
$\alpha,\beta$ be the simple roots of $\Phi$ and denote by $s,r$ the
corresponding simple reflection. We assume that $\alpha$ is a long
root, then $\Phi^+$ consists of $\alpha,\ \beta,\ \alpha+\beta,\
\alpha+2\beta,\ \alpha+3\beta,\ 2\alpha+3\beta.$ We have $W=\{e,\
s,\ r,\ sr,\ rs,\ srs,\ rsr,\ srsr,\ rsrs,\ srsrs,\ rsrsr,\
rsrsrs\}$ and $Y_{\{s\}}=\{e,\ r,\ sr,\ rsr,\ srsr\}$.

Using (2) we have $U_s=U_\alpha$, $U_{sr}=
U_{\alpha+3\beta}U_\beta$, $U_{srs}= U_{2\alpha+3\beta}
U_{\alpha+\beta} U_\alpha$, $U_{srsr}=U_{2\alpha+3\beta}
U_{\alpha+2\beta}U_{\alpha+3\beta}U_\beta$,
$U_{srsrs}=U_{\alpha+3\beta} U_{\alpha+2\beta}U_{2\alpha+3\beta}
U_{\alpha+\beta} U_\alpha$ and denote by $V_i=U_w$ where $l(w)=i$.
By Lemma 2.7 we get
\begin{equation}E_s=\mathfrak KV_1C_s+\mathfrak KV_2rC_s+\mathfrak KV_3srC_s+
\mathfrak KV_4rsrC_s+\mathfrak KV_5srsrC_s,\end{equation} where
$C_s$ is the image of $(1-s)1_{tr}$ in $E_s$.

We show that $E_s$ is generated by any nonzero element in $E_s$ so
that $E_s$ is irreducible. Let $\xi$ be a nonzero element of $E_s$.
Using (23) we have
\begin{equation}\xi=\underset{u\in V_1}{\sum}a_uuC_s+\underset{u\in V_2}{\sum}b_uurC_s+\underset{u\in V_3}{\sum}c_uusrC_s+
\underset{u\in V_4}{\sum}d_uursrC_s+\underset{u\in
V_5}{\sum}e_uusrsrC_s,\end{equation} where $a_u,b_u,c_u,d_u,e_u\in
\mathfrak K$ and only finitely many of them are nonzero.

Clearly there is a positive integer $m$ such that $u\in U_{q^m}$ if
$a_u\ne 0$, $b_u\ne 0$, $c_u\ne 0$, $d_u\ne 0$ or $e_u\ne 0$. We may
require that if $a_u\ne 0$ (resp. $b_u\ne 0$; $c_u\ne 0$; $d_u\ne
0$; $e_u\ne 0$), then $u$ is $\alpha$-regular (resp.
$\beta$-regular; $\alpha$-regular; $\beta$-regular;
$\alpha$-regular).  Otherwise, according to Lemma 2.3, we can choose
some $y \in U$ such that $y\xi$ satisfies the requirement. Then we
consider $y\xi$ instead of $\xi$.

Let $n_s$ (resp. $n_r$) be a representative of $s$ (resp. $r$) in
$N=N_G(T)$. Then $n_s^2$ and $n_r^2$ are in $T$. Since $m$  is
sufficiently large, we may require that  $n_s$ and $n_r$ are in
$G_{q^m}$. Noting that the image of $\eta_S$ where $S=\{s,r\}$ in
$E_s$ is zero, we have
\begin{equation} rsrsrC_s =C_s-rC_s+srC_s-rsrC_s+srsrC_s. \end{equation}

With the assumption of $\xi$, using 2.4(c), 2.4(d), 2.4(e), we get
\medskip

(a) If $a_u\ne 0$ (resp. $b_u\ne 0$; $c_u\ne 0$; $d_u\ne 0$; $e_u\ne
0$), then $n_ruC_s \in UrC_s$ (resp. $n_rurC_s \in UrC_s$;
$n_rusrC_s \in UrsrC_s$; $n_rursrC_s \in UrsrC_s$; $n_rusrsrC_s \in
UrsrsrC_s$).

\medskip

(b) If $a_u\ne 0$ (resp. $b_u\ne 0$; $c_u\ne 0$; $d_u\ne 0$; $e_u\ne
0$), then $n_suC_s= (y-z)C_s$ for some $y,z\in U$ (resp. $n_surC_s
\in UsrC_s$; $n_susrC_s \in UsrC_s$; $n_sursrC_s \in UsrsrC_s$;
$n_susrsrC_s \in UsrsrC_s$).

\medskip

For convenience, set
$$M'=\mathfrak KG\xi,\quad \mathbb{X}=\displaystyle
\sum_{x\in U_{q^m}}x, \quad \phi_h=\displaystyle \underset{u\in
V_i}{\sum}h_u \ \text{for}\  h=a,\ b,\ c,\ d,\ e,\ a',\ b' \dots $$

Let $\xi_1= n_r \xi$ and using (a) $\xi_1$ can be written as
\begin{equation}\xi_1= \underset{u\in V_2}{\sum}b^{(1)}_uurC_s+
\underset{u\in V_4}{\sum}d^{(1)}_uursrC_s+\underset{u\in
U}{\sum}e^{(1)}_uursrsrC_s \in E_s.\end{equation} Here we consider
$rsrsrC_s$ as its image in $E_s$ and for convenience, we do not
write it as $C_s- rC_s+ srC_s-rsrC_s+srsrC_s$. We can also obtain
that $\phi_{b^{(1)}}=\phi_a+\phi_b,\  \phi_{d^{(1)}}=\phi_c+\phi_d,\
\phi_{e^{(1)}}=\phi_{e}$.

We can assume that all $u\in U$ appear in $\xi_1$ are
$\alpha$-regular. Otherwise, we can choose some $y \in U$ such that
$y\xi_1$ satisfies the requirement by Lemma 2.3. Then we consider
$y\xi_1$ instead of $\xi_1$.

With the assumption of $\xi_1$, it is easy to see following result.
\medskip

(c) If $e^{(1)}_u\ne 0$, then $n_s ursrsrC_s= (x-y)rsrsrC_s$ for
some $x,y\in U$.

\medskip

Let $\xi_2 =n_s \xi_1$ and using (b) and (c) we have
\begin{equation}\xi_2 = \underset{u\in V_3}{\sum}c^{(2)}_uusrC_s+
\underset{u\in V_5}{\sum}e^{(2)}_uusrsrC_s+\underset{u\in
U}{\sum}f^{(2)}_{y,z}(y-z)rsrsrC_s, \end{equation} where
$\phi_{c^{(2)}}=\phi_a+\phi_b,\ \phi_{e^{(2)}}=\phi_c+\phi_d$.

We can also assume that  all $u\in U$ appear in $\xi_2$ are
$\beta$-regular. Let $\xi_3=n_r \xi_2$ and using (a), we have
\begin{equation}\xi_3= \underset{u\in V_4}{\sum}c^{(3)}_uursrC_s+
\underset{u\in U}{\sum}e^{(3)}_uursrsrC_s+\underset{y,z\in
U}{\sum}f^{(3)}_{y,z}(y-z)rsrsrC_s,\end{equation} where
$\phi_{c^{(3)}}=\phi_{c^{(2)}}=\phi_a+\phi_b,\
\phi_{e^{(3)}}=\phi_{e^{(2)}}$.

Assume that all $u\in U$ appear in $\xi_3$ are $\alpha$-regular. Let
$\xi_4=n_s\xi_3$ and by (b) and (c) we have
\begin{equation}\xi_4=\underset{u\in V_5}{\sum}c^{(4)}_uusrsrC_s+
\underset{y,z\in U}{\sum}f^{(4)}_{y,z}(y-z)rsrsrC_s
\end{equation} where $\phi_{c^{(4)}}=\phi_{c^{(3)}}=\phi_a+\phi_b$.

Noting $\mathbb{X}u=\mathbb{X}$ for any $u\in U_{q^m}$ we can obtain
$$\mathbb{X} \xi_4 = ( \phi_a+\phi_b)\mathbb{X}srsrC_s \in M'.$$
Thus if $\phi_a+\phi_b \ne 0$ then $srsrC_s\in M'$ by Lemma 2.9 and
we have $M'=E_s.$

Assume that $\phi_a+\phi_b = 0$. Using (27) we can see
$$\mathbb{X} \xi_2 = (\phi_c+\phi_d) \mathbb{X}srsrC_s \in M'.$$
Hence, If $\phi_c+\phi_d \ne 0$, Using Lemma 2.9 we have $srsrC_s\in
M'$ which implies $M'=E_s$.

Now we assume that $\phi_a+\phi_b = 0$ and $\phi_c+\phi_d = 0$.
Using (26) we get
$$\mathbb{X} \xi_1 = \phi_e \mathbb{X} rsrsrC_s \in M'.$$
Hence, if $ \phi_e \ne 0$, then $rsrsrC_s\in M'$ by Lemma 2.9 and we
have $M'=E_s$.

Let $\zeta_1=n_s \xi$ and using (b) we can see that
\begin{equation}\zeta_1= \underset{y,z\in V_1}{\sum}f_{y,z}(y-z)C_s+
\underset{u\in V_3}{\sum}g_uusrC_s+ \underset{u\in
V_5}{\sum}h_uusrsrC_s, \end{equation} where $\phi_g=\phi_b+\phi_c,\
\phi_h=\phi_d+\phi_e.$

We can also assume that  all $u\in U$ appear in $\zeta_1$ are
$\beta$-regular by Lemma 2.3. Let $\zeta_2= n_r \zeta_1$ and using
(a), we have
\begin{equation}\zeta_2= \underset{y,z\in V_2}{\sum}f'_{y,z}(y-z)rC_s+
\underset{u\in V_3}{\sum}g'_uursrC_s+ \underset{u\in
U}{\sum}h'_uursrsrC_s \end{equation} and
$\phi_{g'}=\phi_{g}=\phi_b+\phi_c$.

As before, We can also assume that  all $u\in U$ appear in $\zeta_2$
are $\alpha$-regular by Lemma 2.3. Let $\zeta_3= n_s \zeta_2$ and
using (b) we obtain
\begin{equation}\zeta_3= \underset{y,z\in V_3}{\sum}f''_{y,z}(y-z)srC_s+
\underset{u\in V_5}{\sum}g''_uusrsrC_s+ \underset{y,z\in
U}{\sum}k_{y,z}(y-z)rsrsrC_s,\end{equation} where
$\phi_{g''}=\phi_b+\phi_c$. It is not difficult to see
$$\mathbb{X}\zeta_3 = (\phi_b+\phi_c)\mathbb{X}srsrC_s \in M'.$$ Hence, if
$\phi_b+\phi_c \ne 0$ then $srsrC_s\in M'$ by Lemma 2.9 and we have
$M'=E_s$.

Now we assume that $\phi_b+\phi_c = 0$. Using (30) we have
$$\mathbb{X}\zeta_1= (\phi_d+\phi_e)\mathbb{X}srsrC_s \in M'.$$
Hence, if $\phi_d+\phi_e \ne 0$ then $srsrC_s\in M'$ by Lemma 2.9
and we have $M'=E_s$. By the computation above we find that

\medskip

(d) If one of $\phi_a,\
\phi_b,\ \phi_c,\ \phi_d,\ \phi_e$ is nonzero then $M'=E_s$. Thus we
just need to deal with the case $\phi_a =
\phi_b=\phi_c=\phi_d=\phi_e=0$.

\medskip

Now assume that $a_{u_1}\ne 0$ for some $u_1 \in V_1$. Then by
2.4(c), 2.4(d), 2.4(e), it is easy to see that
$$n_s{u_1}^{-1}\xi= \underset{u\in V_1}{\sum}a^*_uuC_s+\underset{u\in V_2}{\sum}b^*_uurC_s+
\underset{u\in V_3}{\sum}c^*_uusrC_s+ \underset{u\in
V_4}{\sum}d^*_uursrC_s+\underset{u\in V_5}{\sum}e^*_uusrsrC_s$$
where $\phi_{a^*} =-a_{u_1}$. By the result in (d), we have $M'=E_s$
in this case.

Now we assume that all $a_u$ are zero but not all $b_u,c_u,d_u,e_u$
in $\xi$ are zero. Let $\zeta= n_s \xi$ and with the assumption of
$\xi$ we have
\begin{equation}\zeta= \underset{u\in V_3}{\sum}p_uusrC_s +\underset{u\in V_5}{\sum}q_uusrsrC_s,\end{equation}
where $\phi_p = \phi_q =0$ and $p_u$, $q_u$ are not all zero.

We compute the form of $n_sn_rn_susrC_s$ for any $u\in V_3$ and
$n_sn_rn_susrsrC_s$ for any $u \in V_5$. The following results can
be obtained from 2.4(c), 2.4(d), 2.4(e).

\medskip

(e) If $u=1$ then $n_sn_rn_susrC_s= -C_s$. If $u_{\alpha}\ne 1$ then
$n_sn_rn_susrC_s\in UsrsrC_s$. If $u_{\alpha}=1$ but
$u_{\alpha+\beta}\ne 1$ then $n_sn_rn_susrC_s\in UsrC_s$. If
$u_{\alpha}=u_{\alpha+\beta}=1$ but $u_{2\alpha+3\beta}\ne 1$ then
$n_sn_rn_susrC_s= (x-y)C_s$ for some $x,y\in U$.

\medskip

(f) If $u_{\alpha}=u_{\alpha+\beta}=u_{2\alpha+3\beta}=1$ then
$n_sn_rn_susrsrC_s\in UrC_s$. If $u_{\alpha}=u_{\alpha+\beta}=1$ but
$u_{2\alpha+3\beta}\ne 1$ then $n_sn_rn_susrsrC_s\in UsrC_s$. If
$u_{\alpha}=1$ but $u_{\alpha+\beta}\ne 1$ then $n_sn_rn_susrsrC_s
\in UsrsrC_s$. If $u_{\alpha}\ne 1$ then $n_sn_rn_susrsrC_s\in
\mathfrak KUrsrsrC_s$.

\medskip

In the formula (33) of $\zeta$ we assume $p_{u_1}\ne 0$ for some
$u_1$ and we let ${\zeta}^* = n_sn_rn_s{u_1}^{-1}\zeta$,
$${\zeta}^*= \underset{u\in V_1}{\sum}a'_uuC_s+\underset{u\in V_2}{\sum}b'_uurC_s+
\underset{u\in V_3}{\sum}c'_uusrC_s+ \underset{u\in
V_4}{\sum}d'_uursrC_s+\underset{u\in V_5}{\sum}e'_uusrsrC_s.$$ Using
(e), (f) and (25)  we can see  $\phi_{a'}+ \phi_{d'}= -p_{u_1}$
which is nonzero. Then $\phi_{a'}\ne 0$ or $\phi_{d'}\ne 0$. We can
deal with this case by the result in (d) and obtain $M'=E_s$.

The last case is all $p_u$ are zero but  $q_{u_2} \ne 0$ for some
$u_2 \in V_5$ in the formula (33). Let
$$\zeta'= n_s {u_2}^{-1}\zeta= \underset{u\in V_4}{\sum}p'_uursrC_s+\underset{u\in V_5}{\sum}q'_uusrsrC_s$$
and $p'_1= q_{u_2}$ is nonzero.

We compute the form of $\Psi_1=n_sn_rn_sn_rursrC_s$ for any $u\in
V_4$ and $\Psi_2=n_sn_rn_sn_rusrsrC_s$ for any $u \in V_5$. The
following results can be obtained from 2.4(c), 2.4(d), 2.4(e).

\medskip
(g) If $u_{\beta} =u_{\alpha+3\beta}=
u_{\alpha+2\beta}=u_{2\alpha+3\beta}=1$ then $\Psi_1=-C_s$. If
$u_{\beta} =u_{\alpha+3\beta}= u_{\alpha+2\beta}=1$ but
$u_{2\alpha+3\beta}\ne 1$ then $\Psi_1=(x-y)C_s$ for some $x,y\in
U$. For other case, we have $\Psi_1\in \mathfrak KUsrC_s+  \mathfrak
KUsrsrC_s+  \mathfrak KUrsrsrC_s$.

\medskip

(h) We always have $\Psi_2\in \mathfrak KUrsrC_s+  \mathfrak
KUsrsrC_s+  \mathfrak KUrsrsrC_s$.

\medskip

Let $\zeta''= n_sn_rn_sn_r \zeta_1$ and we write $\zeta''$ as
$$\zeta''=\underset{u\in V_1}{\sum}a''_uuC_s+\underset{u\in V_2}{\sum}b''_uurC_s+
\underset{u\in V_3}{\sum}c''_uusrC_s+ \underset{u\in
V_4}{\sum}d''_uursrC_s+\underset{u\in V_5}{\sum}e''_uusrsrC_s.$$
Then by (g), (h) and (25) we can obtain $\phi_{a''}+ \phi_{b''}=
-p'_1 $ which is nonzero. Then $\phi_{a''} \ne 0$ or $\phi_{b''} \ne
0$. We can deal with this case by the result in (d) and obtain
$M'=E_s$.

We can also prove that $E_r$ is irreducible in the similar way. Now
we complete the proof of Theorem 3.1.

\medskip

By Theorem 3.1 the $\mathfrak KG$ module $M(tr)$ has the following
composition series $\displaystyle 0\  \subsetneqq \ St \ \subsetneqq
\ M(tr)_{\{s\}} \ \subsetneqq \ M(tr)_{\{s\}}+M(tr)_{\{r\}} \
\subsetneqq \ M(tr).$

\bigskip

\section{Type $A$ case, I }

In this section  $G$ is assumed to be a
reductive group whose derived group is of type $A_n$. The  Weyl
group then is isomorphic to the symmetric group of $n+1$ letters.
The main result of this section is the following result.

\medskip

{\bf Theorem 4.1}
 Assume that $G$ is a connected reductive group
over $\bar{\mathbb F}_q$ whose derived group is of type $A_n$ and
char\,$\mathfrak K\ne$\,char\,$\bar{\mathbb{F}_q}$. Then

\medskip

(a) the $\mathfrak KG$ module $E_s$ is irreducible for any simple
reflection;

\medskip

(b) the $\mathfrak KG$ modules $E_J$ is irreducible if $J$ is a
maximal proper subset of the $S$. (Recall that $S$ is the set of
simple reflections of the Weyl group of $G$.)

\medskip

{\bf 4.2\ \ } The rest of this section is devoted to prove part (a)
of the theorem above. Part (b) will be proved in next section. We
number the reflections $s_1, s_2, \dots ,s_n$ and the the
corresponding simple roots $\alpha_1,\alpha_2,\dots ,\alpha_n$ as
usual, so the Dynkin diagram of the root system $\Phi$ of $G$ is as
follows:

\medskip

\centerline{\begin{tabular}{lll}
\xymatrix{\underset{1}{\circ}\ar@{-}[r]&\underset{2}{\circ}\ar@{-}[r]
&\underset{3}{\circ}\ar@{-}[r]&\cdots\cdots\ar@{-}[r]&\underset{n-1}{\circ}\ar@{-}[r]&\underset{n}{\circ}}
\end{tabular}}
\def\fK{\mathfrak K}

\medskip

{\bf 4.3\ \ } We first show that $ E_{s_1}$ and $ E_{s_n}$ are simple $\mathfrak KG$-modules.
By symmetry, it is enough to prove that $E_{s_1}$ is irreducible.

Let $C_1$ be the image in  $E=E_{s_1}$ of $(1-s_1)1_{tr}$.  We have
$$Y_{\{s_1\}}=\{e,\ s_2,\ s_3s_2,\dots ,\ s_n\cdots s_3s_2\}.$$ Let $V_i=U_{s_1s_2\cdots s_i}$. Then $V_i=U_{\alpha_1+\alpha_2+\dots +\alpha_i}U_{\alpha_2+\dots
+\alpha_i}\cdots U_{\alpha_i}$. By Lemma 2.7 we get
\begin{align} E&=\mathfrak KUC_1+\mathfrak KUs_2C_1+ \mathfrak KUs_3s_2C_1+ \dots + \mathfrak KUs_n\dots s_3s_2C_1\notag\\
&=\mathfrak KV_1C_1+\mathfrak KV_2s_2C_1+ \mathfrak KV_3s_3s_2C_1+ \dots + \mathfrak KV_ns_n\dots s_3s_2C_1.\end{align}

We show that $E$  is generated by any nonzero element in $E$ so that
$E$  is irreducible. Let $\xi$ be a nonzero element of $E $, Using
formula (34) we see that
\begin{equation}\xi=\underset{u\in V_1}{\sum}a_{1,u}uC_1+ \underset{u\in V_2}{\sum}a_{2,u}us_2C_1+
\dots +\underset{u\in V_n}{\sum}a_{n,u}us_n\dots s_3s_2C_1,\quad
a_{i,u}\in\fK,\end{equation} and only finitely many of the coefficients are
nonzero.

By Lemma 2.3 we may assume that $u$ is $\alpha_1$-regular whenever
$a_{1,u}$ is nonzero. Otherwise we replace $\xi$ by suitable $y\xi.$
Also we may choose $m$ sufficiently large so that all $n_i=n_{s_i}$
are in  $G_{q^m}$. By 2.4(e) we have the following result.

\medskip

(a) If $a_{1,u}\ne 0$, then $n_1uC_1=(v-1)C_1$ for some $v\in V_1$.

\medskip

Noting that the image of $(1-s_1-s_2+s_1s_2+s_2s_1-s_1s_2s_1)1_{tr}$ in $E$ is zero, we get
\begin{align}n_1s_2C_1&=s_1s_2C_1=s_2C_1-C_1,\notag\\
n_1s_i\cdots s_2C_1&=s_i\cdots s_3 s_1 s_2C_1= s_i\cdots
s_2C_1-C_1,\quad\text{if }i\ge 2.\end{align} Using 2.1(a) and
formula (36) one may verify the following claim.

\medskip

(b) If $i\ge 2$, then there exists
$y\in V_{i,q^m}$ such that $$n_1us_i \dots s_2C_1=y(s_i \dots s_2C_1 -
C_1).$$

\medskip

Set
$$M'=\fK G\xi,\quad \mathbb{X}=\underset{x\in U_{q^m}}{\sum}x,\quad \phi_i=\underset{u\in V_i}{\sum}a_{i,u},\quad \phi=\phi_1 + \phi_2+  \dots +\phi_n.$$
Since $\mathbb{X}u=\mathbb{X}$ for any $u\in U_{q^m}$, Using (a) and
(b) we get
\begin{align*}\mathbb{X}\xi - \mathbb{X}n_1\xi&=
\sum_{i=1}^n\phi_i\mathbb{X}s_i\cdots s_2C_1-\sum_{i=2}^n\phi_i\mathbb{X}(s_i\cdots s_2C_1-C_1)\\
&=\phi \mathbb{X}C_1 \in M',\end{align*} here we understand that
$s_i\cdots s_2=e$ for $i=1$. If $\phi \ne 0$ then $\mathbb{X}C_1$ is
in $M'$. By Lemma 2.9 , $C_1\in M' $, therefore $M'=\mathfrak
KG\xi=\mathfrak KGC_1=E $.

Assume that $\phi= 0$ but $a_{1,u_1} \ne 0$ for some $u_1 \in V_1$. Note that $n_1C_1=-C_1$. Using (a) and (b) we get
$$\mathbb{X}\xi -\mathbb{X}n_1{u_1}^{-1}\xi= a_{1,u_1} \mathbb{X} C_1\in M'.$$
Again, using Lemma 2.9 we see that $M'=\mathfrak KG\xi=\mathfrak
KGC_1=E$ in this case. We have proved the following result.

(c) Let $\xi$ be a nonzero element in $E$ of the form (35) such that
$u\in V_1$ is $\alpha_1$-regular if  $a_{1,u}\ne 0$. Provided that
$\phi\ne 0$ or $a_{1,u}\ne 0$ for some $u\in V_1$, then $\mathfrak
KG\xi=E_s.$

\medskip

Let $i,j$ be integers such that $2\le i<j\le n$. Then $n_iV_{j,q^m}n_i^{-1}=V_{j,q^m}$
and $s_{i+1}C_1=C_1$. So  we have the following result

\medskip

(d) If $2\le i<j\le n$ and $u\in V_{j,q^m}$, then
$$n_ius_j \cdots
s_2C_1=u' s_j \dots s_2C_1$$ for some $u' \in V_{j,q^m}$.

\medskip

Now assume that all $a_{i',u}=0$ whenever $1\le i'<i$ and some
$a_{i,u_i}\ne 0$. By (d), 2.1(a), 2.4 (b), 2.4(c) and 2.4(d) we get
$$\zeta=n_1n_2\cdots n_iu_i^{-1}\xi=\underset{u\in V_1}{\sum}a'_{1,u}uC_1+
\underset{u\in V_2}{\sum}a'_{2,u}us_2C_1+ \dots +\underset{u\in
V_n}{\sum}a'_{n,u}us_ns_{n-1}\dots s_2C_1,$$ where all $a'_{i,u}$ are in $\fk$ and
$a'_{1,1}=a_{i,u_i}\ne 0$.

By Lemma 2.3, we may choose $y\in V_1$ such that $yu$ is
$\alpha_1$-regular when $u\in V_1$ and $a'_{1,u}\ne0$. By (c) we
have
 $\mathfrak KGy\zeta=E.$ Therefore $\mathfrak KG\xi=E$ in this
case.

We have proved that $E_{s_1}$ is irreducible. By symmetry, $E_{s_n}$ is also irreducible.

\medskip

{\bf 4.4}\ \  Now we fix an integer $i'$ such that $2\le i' \le n-1$
and set $s=s_{i'}$. In this subsection we show that $E_s$ is
irreducible.

Using Lemma 2.7 to describe $E_s$ it suffices to describe $Y=
Y_{\{s\}}$. We set $m=i'-1$ and $l=n-i'$. For convenience, denote by
$r_j=s_{i'-j}$ for $j=1,2,\dots,m$, $t_k=s_{i'+k}$ for
$k=1,2,\dots,l$. The Dynkin diagram of $W$ is as follows:

\bigskip

\centerline{\begin{tabular}{lll}
&\xymatrix{\circ_{r_m}\ar@{-}[r]&\circ\cdots\cdots\circ_{r_2}\ar@{-}[r]&\circ_{r_1}\ar@{-}[r]&\circ_{s}\ar@{-}[r]&\circ_{t_1}\ar@{-}[r]&\circ\cdots\cdots\circ_{t_{l-1}}\ar@{-}[r]&\circ_{t_l}}
\end{tabular}}

\bigskip

{\bf Lemma 4.5\ \ } Keep the notations above. Denote by
$\sigma_i=r_ir_{i-1}\dots r_2r_1$, $i=1,2,\dots,m$ and $\tau_j=t_j
t_{j-1}\dots t_2t_1$, $j=1,2,\dots,l$. (Convention: $\sigma_0$ and
$\tau_0$ are the neutral element $e$ in $W$.) Then $Y$ consists of
following elements
$$\tau_{j_k}\sigma_{i_k}s\tau_{j_{k-1}}\sigma_{i_{k-1}}\dots s\tau_{j_2}\sigma_{i_2}s\tau_{j_1}\sigma_{i_1}$$
where $m\ge i_1>i_2>\dots>i_k \geq 0, l\ge j_1>j_2>\dots>j_k \geq
0.$

\medskip
Proof. Denote by $Z$ the element listed in the lemma. We just need
to prove $Y=Z$.

As the notation in 2.5, $X_J=\{x\in W\mid x \ \text{has minimal
length in}\  xW_J\}$ for a subset $J$ of $S$. In this proof of the
lemma, we always set $J=S\setminus\{s\}$. Using  2.1.1 proposition
in [GP], we have $X_J= Ys\cup\{e\}$ and $X_J$ is a left coset
representatives of $W_J$ in $W$.

Without lost of generality, we assume $m\le l$. We denote by
$\delta_J$ the following element
$$\delta_J=\tau_{l-m}s\tau_{l-m+1}\sigma_1\dots s\tau_{l-1}\sigma_{m-1}s\tau_l\sigma_m,$$
which is the longest element of $Z$. We show $\delta_J$ is also the
longest element of $Y$. By easy computation, we know
$l(\delta_J)=ml+m+l$. On the other hand, we can get
$$l(w_0)-l(w_J)=\frac{n(n+1)}{2}-(\frac{m(m+1)}{2}+\frac{l(l+1)}{2})= ml+m+l+1.$$
It suffice to prove $\delta_J$ is in $Y$, then we can see
$\delta_Jsw_J=w_0$ and $\delta_J$ is the longest element of $Y$.

We regard the simple reflection $s_j$ as exchanging the two letters
$j$ and $j+1$. Thus $\delta_J$ in $Y$ means $\delta_Js$ in $X_J$
which is equivalent to the following two conditions:

\medskip
(a) $\delta_Js(i')> \delta_Js(i'+1)$.

\medskip

(b) For $k\ne i'$, $\delta_Js(k)<\delta_Js(k+1)$.

\medskip
Recall $i'$ is a fixed integer and $l=n-i'$ as we mentioned before.
By careful computation, we have the following result.
\medskip

(c) For $k\le i'$, $\delta_Js(k)=l+k+1$.

\medskip

(d) For $k\ge i'$, $\delta_Js(k)=k- i'$.

\medskip

Then (c) and (d) imply $\delta_Js$ satisfy the two condition (a) and (b).
Therefore $\delta_J$ is the unique longest element in $Y$.

It is not difficult to see that for any $w\in Z$,
$l(\delta_Jw^{-1})=l(\delta_J)-l(w)$. Using the 2.2.1 Lemma in [GP],
we can get $Z\subseteq Y$.(Here we use left coset representative
while [GP] used right coset representative.) On the other hand, by
the Algorithm C (p.46) in [GP], we can see  $Z$ is exactly the set
$Y$.

The lemma is proved.\qed

\bigskip

We will show that $E_s$  is generated by any nonzero element in
$E_s$ so that $E_s$  is irreducible. Let $\xi$ be a nonzero element
of $E_s$. By Lemma 2.7 and Lemma 4.5 we have
\begin{equation}\xi=\underset{w\in Y}{\sum} \underset{x\in
V_w}{\sum}a_{w,x}xwC_s,\quad a_{w,x}\in\mathfrak K,\end{equation}
where $V_w= U_{sw^{-1}}$ and only finitely many of the coefficients
are nonzero.

\medskip
By Lemma 2.3 we can assume that $x$ appears in $\xi$ is
$\Delta_{sw^{-1}}$-regular wherever $a_{w,x}$ is nonzero. Otherwise,
we replace $\xi$ by suitable $y\xi$. This process is called
regularlization to $\xi$ and we denote $R(\xi)$ for one element
$y\xi=\underset{w\in Y}{\sum} \underset{z\in
V_w}{\sum}a'_{w,z}zwC_s$ such that $z$ appears in $y\xi$ is
$\Delta_{sw^{-1}}$-regular wherever $a'_{w,z}$ is nonzero.

\medskip

We may choose $m$ sufficiently large so that $n_j=n_{s_j}$ are in
$G_{q^m}$. In the following, we write $s_jxwC_s$ instead of
$n_jxwC_s$ for convenience. Given simple reflections $s_1,s_2,\dots,
s_k$, we denote by
$$R(n_1n_2\dots n_k\xi)=R(n_1(\dots R(n_{k-1}R(n_k\xi)))).$$ We can
also just denote it by $R(s_1s_2\dots s_k\xi)$.

\medskip

For convenience, we classify the elements in $Y$ by the different
endings of the elements in $Y$. For $k$ a positive integer we denote
by
\medskip
$$A_k=\{\sigma_{i_k}s\tau_{j_{k-1}}\sigma_{i_{k-1}}\dots s\tau_{j_1}\sigma_{i_1}\mid m\ge i_1>\dots>i_k >0, l\ge j_1>\dots>j_{k-1} > 0\}.$$
$$B_k= \{\tau_{j_k}s\tau_{j_{k-1}}\sigma_{i_{k-1}}\dots s\tau_{j_1}\sigma_{i_1}\mid m\ge i_1>\dots>i_{k-1}>0,
l\ge j_1>\dots>j_k > 0\}.$$
$$C_k= \{\tau_{j_k}\sigma_{i_k} s\tau_{j_{k-1}}\sigma_{i_{k-1}}\dots s\tau_{j_1}\sigma_{i_1} \mid m\ge i_1>\dots>i_k>0,
l\ge j_1> \dots>j_k>0\}.$$
$$D_{k+1}=\{s\tau_{j_k}\sigma_{i_k}\dots s\tau_{j_1}\sigma_{i_1}\mid m\ge i_1>\dots>i_k>0,
l\ge j_1>\dots>j_k>0\}.$$ and $D_1$ consist of the neutral element $e$.

\medskip

For convenience, we set
$$M'=\fK G\xi,\quad  \mathbb{X}=\underset{x\in U_{q^m}}{\sum}x,\quad \phi_w= \underset{x\in
V_w}{\sum}a_{w,x}.$$ We also denote by $\phi_{H_k}=\underset{w\in
H_k}{\sum}\phi_w$ for $H=A,\ B,\ C,\ D $ and $\phi_k=
\phi_{A_k}+\phi_{B_k}+\phi_{C_k}+\phi_{D_k}$.

\medskip

{\bf(Step 1)\ \ } With the assumption of $\xi$ and by 2.4 (c), 2.4
(d), 2.4 (e),  for $w\in Y, x\in V_w$ we have the following result.

\medskip

(a)If $s_jw> w$, then $s_jxwC_s=x's_jwC_s$ for some $x'\in U$.

\medskip

(b)If $s_jw< w$, then $s_jxwC_s=x''wC_s$ for some $x''\in U$.

\medskip
(c)For $s_j\ne s, r_1, t_1$,  we get $s_jC_s=C_s$.
\medskip

 Denote by $w_a$ for the following element
$$(r_m\dots r_1)(r_m\dots r_2)\dots (r_mr_{m-1})r_m(t_l\dots
t_1)(t_l\dots t_2)\dots (t_lt_{l-1})t_l.$$

We consider the element $\zeta=R(w_a\xi)$. We set $$w_c=(t_l\dots t_2
t_1 r_m\dots r_2r_1)s.$$ Using (a), (b), (c) and the exchange
relations of Weyl group, we can see $\zeta$ has the following form
\begin{equation}\zeta= \underset{k\le min\{m,l\}}{\sum}
\underset{x\in U }{\sum}a_{k,u}x{w_c}^k(t_l\dots t_2t_1 r_m\dots
r_2r_1)C_s.\end{equation}

The elements $(w_c)^k(t_l\dots t_1 r_m\dots r_1)$ may not be in $Y$.
However we can express $\zeta$ of this form since we have the
following lemma.

\medskip

{\bf Lemma 4.6\ } The elements ${w_c}^{k+1}$ is a reduced expression
for $k\le \min\{m,l\}$.

\medskip
Proof. We just need to prove ${w_c}^k$ is a reduced expression when
$k=\min \{m,l\}+1$. For the convenience of computation, we go back
to the notation of $\{s_j\}$. Without lost of generality, we assume
$m\le l$ which implies $n \geq 2m+1$. Let $w={w_c}^m=((s_n\dots
s_{m+2})(s_1\dots s_{m})s_{m+1})^{m+1}$ and we want to prove $w$ is
a reduced expression. Since the length of $w$ equals to the number
of positive roots transformed by $w$ into negative roots. We use
this to prove the lemma by computation.

Let $s_i= (i,i+1)$ as exchanging two letters $i$ and $i+1$. The
positive roots can be written as $\{\varepsilon_i-\varepsilon_j\mid
i<j \}$. Firstly, by computation we can see that

\begin{equation}w_c(\varepsilon_i)= \left\{\begin{aligned}
\varepsilon_{i+1} \ \ \ \ \ \ \ \ \ \ \ \ \ \text{if} \ 1\le i \le m.  \\
\varepsilon_{n+1} \ \ \ \ \ \ \ \ \ \ \ \ \ \text{if} \  i=m+1. \\
\varepsilon_1 \ \ \ \ \ \ \ \ \ \ \ \ \ \ \ \text{if} \  i=m+2. \\
\varepsilon_{i-1} \ \ \text{if} \ m+3 \le i \leq n+1.\\
\end{aligned} \right.
\end{equation}

So it is not hard to compute $w(\varepsilon_i)$,

\begin{equation}w(\varepsilon_i)= \left\{\begin{aligned}
\varepsilon_{n+2-i} \ \ \ \ \ \ \ \ \ \ \ \ \ \ \ \ \ \ \ \ \text{if} \ 1 \le i \le m+1. \\
\varepsilon_{2m+3-i} \ \ \ \ \ \ \ \ \ \ \  \text{if} \ m+2 \le i \le 2m+2. \\
\varepsilon_{i-(m+1)} \ \ \ \ \ \ \ \ \ \ \  \text{if} \ 2m+3 \le i \le n+1. \\
\end{aligned} \right.
\end{equation}

\medskip

By above computation, when we compute the number of  positive roots
transformed by $w$ into a negative roots ,we need to check it case
by case. We consider a positive root $\varepsilon_i -\varepsilon_j$
for $i < j$ and compute the number of pairs $(i,j)$ such that
$w(\varepsilon_i -\varepsilon_j)$ is a negative root.

\medskip

(a) When $1 \leq i \leq m+1$.

When $1 \leq j \leq m+1$, then $i\leq j$ and $n+2-i \geq n+2-j $.
The number of these pairs of $(i,j)$ is $m(m+1)/2$.

When $ m+2 \leq j \leq 2m+2$, then $i\leq j$ and $n+2-i \geq
2m+3-j$. The number of these pairs of $(i,j)$ is $(m+1)^2$.

When $2m+3 \leq j \leq n+1$, then $i\leq j$ and $n+2-i \geq
j-(m+1)$. The number of these pairs of $(i,j)$ is $(n-2m-1)(m+1)$.

\medskip

(b) When $m+2 \leq i \leq 2m+2$.

When $m+2 \leq j \leq 2m+2$, then $i\leq j$ and $2m+3-i \geq 2m+3-j
$. The number of these pairs of $(i,j)$ is $m(m+1)/2$.

When $2m+3 \leq j \leq n+1$, then $i\leq j$ and $2m+3-i \geq
j-(m+1)$. Hence $j < 3m+2-i$ and it is a contradiction.There is no
such pair of $(i,j)$.

\medskip

(c) When $i\geq 2m+3 $, then $j \geq 2m+3$. There is no such pair of
$(i,j)$.

\medskip

Combining (a), (b) and (c), the number of positive roots transformed
by $w$ into negative roots is $$\frac{m(m+1)}{2}+ (m+1)^2+
(n-2m-1)(m+1)+ \frac{m(m+1)}{2}= n(m+1).$$ This number is also the
length of $w$ and the lemma is proved. \qed

\medskip

{\bf(Step 2)\ \ }Keep the notation of  $w_c$. The element $\zeta$ in
Step 1 is of the formula (38).

We consider the element $R({w_c}^h\zeta)$ for some integer $h$. By
2.4 (e),  if we use $\mathbb{X}$ multiply on $R({w_c}^h\zeta)$, we
can see that the bigger $h$ is, the more elements in $\zeta$ will be
killed. This phenomenon can be easy to see in the case of type $G_2$
in subsection 3.5. Moreover we can consider the element $R((t_i\dots
t_2t_1 r_j\dots r_2r_1s){w_c}^h\zeta)$ for some integers $i,j,h$. By
the same philosophy we can see the following result.

\medskip

If $\phi_k$ is nonzero for some integer $k$, then $C_s\in M'$ which
implies $\mathfrak KG\xi =E_s$.

\medskip

In the following, we assume $\phi_k$ is zero for all $k\ge 1$.

\medskip

{\bf(Step 3)\ \ } In this step, when we consider $xwC_s$ we assume
$x\in U$ is an $\Delta_{sw^{-1}}$-regular element. Using the result
in subsection 2.4, the following result is easy to see. The
notations $\sigma_i=r_i\dots r_2r_1$ and $\tau_j=t_j\dots t_2t_1$
are as before.

\medskip

(a) We have $suC_s=(x-1)C_S$ for some $x\in U$. $su\sigma_iC_s=
y(\sigma_iC_s-C_s)$ for some $y \in U$. $su\tau_jC_s=
z(\tau_jC_s-C_s)$ for some $z \in U$.

\medskip
(b) We compute $suA_kC_s$ for $k\ge 2$ and we have
$$su\sigma_{i_k}s\tau_{j_{k-1}}\sigma_{i_{k-1}}\dots s\tau_{j_1}\sigma_{i_1}C_s
=u'\sigma_{i_k}s\tau_{j_{k-1}}\sigma_{i_{k-1}}\dots
s\tau_{j_1}\sigma_{i_1}C_s$$ for some $u'\in U$.

\medskip
(c) We compute $suB_kC_s$ for $k\ge 2$ and we have
$$su\tau_{j_k}s\tau_{j_{k-1}}\sigma_{i_{k-1}}\dots s\tau_{j_1}\sigma_{i_1}C_s
=u''\tau_{j_k}s\tau_{j_{k-1}}\sigma_{i_{k-1}}\dots
s\tau_{j_1}\sigma_{i_1}C_s$$ for some $u''\in U$.

\medskip

(d) We compute $suC_kC_s$ for $k\ge 2$ and we have
$$su\tau_{j_k}\sigma_{i_k} s\tau_{j_{k-1}}\sigma_{i_{k-1}}\dots s\tau_{j_1}\sigma_{i_1} C_s
=u^*s\tau_{j_k}\sigma_{i_k} s\tau_{j_{k-1}}\sigma_{i_{k-1}}\dots
s\tau_{j_1}\sigma_{i_1} C_s$$ for some $u^*\in U$.

\medskip

(e) We compute $suD_kC_s$ for $k\ge 2$ and we have
$$sus\tau_{j_k}\sigma_{i_k}\dots s\tau_{j_1}\sigma_{i_1}C_s
 =u^{**}s\tau_{j_k}\sigma_{i_k}\dots s\tau_{j_1}\sigma_{i_1}C_s$$ for some  $u^{**}\in U$.

\medskip
As the notation before, $m,l$ is the  number of reflections
$\{r_i\}$ and $\{t_j\}$. If $m=l$ then there exists a maximal
integer $h$ such that $D_h$ is nonempty but $A_h,B_h,C_h$ are all
empty. If $m>l$ then there exists a maximal integer $h$ such that
$A_h,D_h$ are nonempty but $B_h,C_h$ are both empty. If $m < l$ then
there exists a maximal integer $h$ such that $B_h,D_h$ are nonempty
but $A_h,C_h$ are both empty.

In each case above, we fix such an integer $h$ and consider $s\xi$.
Using (a), (b), (c), (d), (e) we can see
$$suA_kC_s\subseteq UA_kC_s,\ \
suB_kC_s\subseteq UB_kC_s, \ \ suD_kC_s\subseteq UD_kC_s \ \  \text{for}\  k\ge
2$$ while $suC_kC_s\subseteq UD_{k+1}C_s$ for $k\ge 1$. Using the result
in step 2 that $\phi_k$  nonzero implies $M' =E_s$, then we see that
if $\phi_{C_h}$ is nonzero we have $M' =E_s$. Moreover by the same
reason we can see the following result.

(f) If $\phi_{C_k}$ is nonzero for any $k\ge 1$ we  have $M' =E_s$.

In the following, we assume $\phi_{C_k}$ and $\phi_k$ are zero.
Otherwise we have $M'=E_s$.

\medskip
{\bf(Step 4)\ \ } Assume that $a_{e,x_0}\ne 0$ of $\displaystyle
\xi=\underset{w\in Y}{\sum} \underset{x\in V_w}{\sum}a_{w,x}xwC_s$
for some $x_0\in U_{\alpha_s}$.
  Using 2.1(a) and 2.4 (d), we observe that if $x_\alpha =1$ then $sxD_kC_s \subset UC_{k-1}C_s$.
Otherwise, $sxD_kC_s \subset UD_kC_s$, where $k\ge 2$. Let
$$s{x_0}^{-1}\xi= \underset{w\in
Y}{\sum} \underset{x\in V_w}{\sum}a'_{w,x}xwC_s.$$

Using (a) in step (3), we obtain $\phi_{A_1}+ \phi_{B_1}+\phi_{D_1}=
-a_{e, x_0}$ which is nonzero. Using the result in step 2 and (f) in
step 3 we can get that $M' =E_s$.

\medskip
{\bf(Step 5)\ \ } In this step, we compute the form of $t_1xwC_s$
without any assumption of $x\in U$. Using 2.4 (c), 2.4 (d) and 2.4
(e) we can have the following results.

\medskip

(a) We compute the form of $t_1uA_kC_s$. In following equations, $u'$
is some element in $U$.
\medskip

(a1) When $k=1$, then $t_1ur_i\dots r_1C_s= u't_1 r_i\dots r_1C_s$.

\medskip

(a2) When $k= 2$ we have following two cases.

If $j_1\ge 2$ then  $t_1u\sigma_{i_2} s \tau_{j_{1}}\sigma_{i_1}C_s
= u't_1\sigma_{i_2} s \tau_{j_{1}}\sigma_{i_1}C_s.$

If $j_1 = 1$ then $t_1u\sigma_{i_2} s t_1 \sigma_{i_1}C_s =u'
(\sigma_{i_2} s t_1 \sigma_{i_1}C_s - t_1\sigma_{i_2}C_s +
\sigma_{i_2}C_s)$

\medskip
(a3) When $k\ge 3$ we also have following two cases.

If $j_{k-1}\ge 2$ then
$$t_1u\sigma_{i_k}s\tau_{j_{k-1}}\sigma_{i_{k-1}}\dots s\tau_{j_1}\sigma_{i_1}C_s
=u't_1\sigma_{i_k}s\tau_{j_{k-1}}\sigma_{i_{k-1}}\dots
s\tau_{j_1}\sigma_{i_1}C_s.$$

If $j_{k-1}=1$ then
$$t_1u\sigma_{i_k}st_1\sigma_{i_{k-1}}\dots s\tau_{j_1}\sigma_{i_1}C_s
=u'\sigma_{i_k}st_1\sigma_{i_{k-1}}\dots
s\tau_{j_1}\sigma_{i_1}C_s.$$

\medskip
(b) We compute the form of $t_1uB_kC_s$.

\medskip

(b1) If $j_k \ge 2$ or $u_{\alpha_{t_1}}\ne 1$ then we have
$$t_1u\tau_{j_k}s\tau_{j_{k-1}}\sigma_{i_{k-1}}\dots s\tau_{j_1}\sigma_{i_1}C_s
=u'\tau_{j_k}s\tau_{j_{k-1}}\sigma_{i_{k-1}}\dots
s\tau_{j_1}\sigma_{i_1}C_s$$ for some $u'\in U$.

\medskip

(b2) If $j_k = 1$ and $u_{\alpha_{t_1}}= 1$ then
$$t_1ut_1 s\tau_{j_{k-1}}\sigma_{i_{k-1}}\dots s\tau_{j_1}\sigma_{i_1} C_s
=u''s\tau_{j_{k-1}}\sigma_{i_{k-1}}\dots s\tau_{j_1}\sigma_{i_1}
C_s$$ for some  $u''\in U$

\medskip
(c) We compute the form of $t_1uC_kC_s$.

\medskip

(c1) If $j_k \ge 2$ or $u_{\alpha_{t_1}}\ne 1$we have
$$t_1u\tau_{j_k}\sigma_{i_k} s\tau_{j_{k-1}}\sigma_{i_{k-1}}\dots s\tau_{j_1}\sigma_{i_1} C_s
=u'\tau_{j_k}\sigma_{i_k} s\tau_{j_{k-1}}\sigma_{i_{k-1}}\dots
s\tau_{j_1}\sigma_{i_1} C_s$$ for some $u'\in U$.

\medskip

(c2) If $j_k= 1$ and $u_{\alpha_{t_1}}= 1$ then
$$t_1ut_1\sigma_{i_k} s\tau_{j_{k-1}}\sigma_{i_{k-1}}\dots s\tau_{j_1}\sigma_{i_1}C_s
=u''\sigma_{i_k} s\tau_{j_{k-1}}\sigma_{i_{k-1}}\dots
s\tau_{j_1}\sigma_{i_1}C_s$$ for some $u''\in U$.

\medskip
(d) We compute the form of $t_1uD_kC_s$. The result is very similar to
(a) and we assume that $u^*$ is some element in $U$.

\medskip

(d1) When $k=1$, then $t_1uC_s= u^*t_1C_s$.

\medskip

(d2) When $k= 2$ we have following two cases.

If $j_1\ge 2$ then  $t_1us \tau_{j_1}\sigma_{i_1}C_s=u^*t_1 s
\tau_{j_1}\sigma_{i_1}C_s.$

If $j_1 = 1$ then $t_1u s t_1 \sigma_{i_1}C_s =u^* (s t_1
\sigma_{i_1}C_s- t_1 C_s + C_s)$.

\medskip

(d3) When $k\ge 3$ we also have following two cases:

If $j_k\ge 2$ then
$$t_1us\tau_{j_k}\sigma_{i_k}\dots s\tau_{j_1}\sigma_{i_1}C_s =u^*
t_1s\tau_{j_k}\sigma_{i_k}\dots s\tau_{j_1}\sigma_{i_1}C_s.$$

If $j_k=1$ then
$$t_1u s\tau_{j_k}\sigma_{i_k}\dots s\tau_{j_1}\sigma_{i_1}C_s
=u^*s\tau_{j_k}\sigma_{i_k}\dots s\tau_{j_1}\sigma_{i_1}C_s.$$

\medskip

{\bf(Step 6)\ \ } In this step, we compute the form of $t_ixwC_s$
for $i\ge 2$ without any assumption of $x\in U$. It is easy to note
that we just need to compute $t_iuB_kC_s$. The form of of $t_iuA_kC_s,
t_iuC_kC_s,t_iuD_kC_s$ can be obtained from $t_iuB_kC_s$.

When we consider the following form $$t_iu
\tau_{j_k}s\tau_{j_{k-1}}\sigma_{i_{k-1}}\dots
s\tau_{j_1}\sigma_{i_1}C_s$$ we have the following results by 2.4
(c), 2.4 (d), 2.4 (e).

\medskip

(a) If $i> j_h+1$ for all $h=1,2,\dots,k$ then
$$t_iu \tau_{j_k}s\tau_{j_{k-1}}\sigma_{i_{k-1}}\dots s\tau_{j_1}\sigma_{i_1}C_s =
u'\tau_{j_k}s\tau_{j_{k-1}}\sigma_{i_{k-1}}\dots
s\tau_{j_1}\sigma_{i_1}C_s$$ for some element $u'\in U$.

Otherwise, there is an integer $a$ such that $j_{a+1}+1<i\le j_a+1$.
Without lot of generality, we can assume $i\le j_k+1$.

\medskip
(b) If $i= j_k+1$ and $j_{k-1}> j_k+1$ (including the case $k=1$)
then
$$t_iu \tau_{j_k}s\tau_{j_{k-1}}\sigma_{i_{k-1}}\dots s\tau_{j_1}\sigma_{i_1}C_s =
u''\tau_{j_k+1}s\tau_{j_{k-1}}\sigma_{i_{k-1}}\dots
s\tau_{j_1}\sigma_{i_1}C_s$$ for some element $u''\in U$.

\medskip

(c) If $i= j_k+1$ and $j_{k-1}= j_k+1$ then we have two cases:

\medskip

(c1) When $k\ge 3$ then
$$t_iu \tau_{j_k}s\tau_{j_{k-1}}\sigma_{i_{k-1}}\dots s\tau_{j_1}\sigma_{i_1}C_s
=x\tau_{j_k}s\tau_{j_{k-1}}\sigma_{i_{k-1}}\dots
s\tau_{j_1}\sigma_{i_1}C_s$$ for some element $x \in U$.

\medskip

(c2) When $k= 2$ then
$$t_{{j_2}+1}u \tau_{j_2}s\tau_{j_1}\sigma_{i_1} C_s =x'(\tau_{j_2}s\tau_{j_1}\sigma_{i_1} C_s-
\tau_{j_1}C_s+ \tau_{j_2}C_s)$$ for some element $x' \in U$.

\medskip

(d) If $i= j_k$ then we also have two cases:

\medskip

(d1) When $u_{\alpha_{t_{j_k}}}\ne 1$ then
$$t_iu
\tau_{j_k}s\tau_{j_{k-1}}\sigma_{i_{k-1}}\dots
s\tau_{j_1}\sigma_{i_1} C_s = y
\tau_{j_k}s\tau_{j_{k-1}}\sigma_{i_{k-1}}\dots
s\tau_{j_1}\sigma_{i_1}C_s$$ for some element $y \in U$.

\medskip

(d2) When $u_{\alpha_{t_{j_k}}}= 1$ then
$$t_iu\tau_{j_k}s\tau_{j_{k-1}}\sigma_{i_{k-1}}\dots s\tau_{j_1}\sigma_{i_1}C_s = y'
\tau_{j_k-1}s\tau_{j_{k-1}}\sigma_{i_{k-1}}\dots
s\tau_{j_1}\sigma_{i_1}C_s$$ for some element $y' \in U$.

\medskip
(e) If $i< j_k$ then
$$t_iu \tau_{j_k}s\tau_{j_{k-1}}\sigma_{i_{k-1}}\dots s\tau_{j_1}\sigma_{i_1}C_s = z
\tau_{j_k}s\tau_{j_{k-1}}\sigma_{i_{k-1}}\dots
s\tau_{j_1}\sigma_{i_1}C_s$$ for some element $z \in U$.

\medskip

{\bf(Step 7)\ \ }We can also compute the form of
$r_iuA_kC_s$, $r_iuB_kC_s$, $r_iuC_kC_s$, $r_iuD_kC_s$ which is similar to the result in
(Step 5) and (Step 6). These computations seem to be complicated.
However using these results in step 5 and step 6, we can have the
following conclusion which is very important to our proof later.

Suppose $t$ is a reflection in $\{r_1,r_2,\dots,r_m, t_1,t_2,\dots,
t_l\}$. Let $w\in Y$ and $x$ be a general element in $U$.

\medskip
 If $l(tw)=l(t)+l(w)$ then $txwC_s$ can only have three cases.

\medskip
(a) We have $txwC_s = utwC_s$ for some element $u \in U$, where $tw\in Y$.

\medskip
(b) We have $txwC_s = y w C_s$ for some element $y\in U$.

\medskip
(c) We have $txwC_s = z (w C_s -w_1C_s+w_2C_s)$ for some element
$z\in U$, $w_1,w_2\in W$ where $w_1,w_2\in Y$ satisfy $l(w_1),l(w_2)<
l(w)$.

\medskip
If $l(tw)=l(w)-l(t)$ then $txwC_s$ can only have two cases.

\medskip
(d) When $x_{\alpha_t}=1$ we have $txwC_s = x' twC_s$ for some
$x'\in U$.

\medskip

(e) When $x_{\alpha_t}\ne 1$ we have $txwC_s = x''wC_s$ for some
$x''\in U$.

\medskip

{\bf(Step 8)\ \ }As many things we prepare in the steps before we
can prove part (a) of Theorem 4.1. Let $\xi$ be a nonzero element in
$E_s$ of the formula (37). Using the observation in (Step 4), we
just need to prove that if there exists a element $w_1$ such that some $a_{{w_1},x}$ are nonzero and
$l(w_1)=d\ge 1$ is least, then we can construct an new element
$$\displaystyle  \xi'=
\underset{w\in Y}{\sum} \underset{x\in V_w}{\sum}a'_{w,x}xwC_s\quad
\in \mathfrak KG\xi$$ such that some $a'_{{w_2},x}$ are nonzero
where $l(w_2)<d$.

\medskip

Using result in  (Step 3) and (Step7), this is not difficult to
prove by induction. As in the assumption,  let $\xi$ be a nonzero
element of the formula (37) and we can assume that $x$ appears in
$\xi$ is $\Delta_{sw^{-1}}$-regular wherever $a_{w,x}$ is nonzero as
before. According to the inductive hypothesis, we assume all
$a_{w,x}=0$ for all $l(w)<d$ but for some $w_1$ and $l(w_1)=d$ such
that $a_{{w_1},{x_1}}$ is nonzero.

Let $t$ be a reflection in
$\{s,r_1,r_2,\dots,r_m,t_1,t_2,\dots,t_l\}$ such that
$l(tw_1)=l(w_1)-l(t)$ and $tw_1 \in Y$. We use the same method in
subsection 3.4. By the results in step 3 and step7 we know that if
$t{x_1}^{-1}\xi$ satisfies our requirement then we can let $\xi' =
t{x_1}^{-1}\xi$. Otherwise, there exists one element $y \in
U_{\alpha_t}$ such that $y x_{\alpha_t}\ne 1$ for $x\in V_{w_1}$,
$a_{w_1,x}\ne 0$ and moreover, this element can make $ty\xi$
satisfy our requirement.

Since $\xi' \in \mathfrak KG\xi$, by induction if we prove $E_s=\mathfrak KG\xi'$ then
$E_s=\mathfrak KG\xi$.

\medskip
Now we complete the proof and obtain $E_s$ is irreducible for any
simple reflections.\qed

\bigskip

\section{type $A$ case, II }
In this section we give the proof of part (b) in Theorem 4.1. We
prove that $\mathfrak KG$ modules $E_J$ is irreducible if $J$ is a
maximal proper subset of the $S$ (the set of simple reflections of
the Weyl group of $G$).

\medskip

{\bf 5.1}\ \ As in the subsection 4.2, we number the reflections
$s_1,s_2,\dots, s_n$  and the the corresponding simple roots
$\alpha_1, \alpha_2,\dots, \alpha_n$ as usual. The Dynkin diagram of
the root system $\Phi$ of $G$ is as follows:

\medskip

\centerline{\begin{tabular}{lll}
\xymatrix{\underset{1}{\circ}\ar@{-}[r]&\underset{2}{\circ}\ar@{-}[r]
&\underset{3}{\circ}\ar@{-}[r]&\cdots\cdots\ar@{-}[r]&\underset{n-1}{\circ}\ar@{-}[r]&\underset{n}{\circ}}
\end{tabular}}
\def\fK{\mathfrak K}

\medskip

Firstly we consider the special case $J=\{2,3,\dots,n \}$ and prove
$E_J$ is irreducible. We show that $E_J$  is generated by any
nonzero element in $E_J$ so that $E_J$ is irreducible. Let $C_J$ be
the image of $\eta_J$ in $E_J$. We have
$$Y_J=\{e,\ s_1,\ s_2s_1,\
\dots,\ s_{n-1}\dots s_2s_1\}.$$ Denote by $V_i=U_{w_Js_1\dots
s_{i-1}}$. By Lemma 2.7 we get
\begin{equation}E_J=\fK V_1C_J+\fK V_2s_1C_J+ \fK V_3s_2s_1C_J+ \dots + \fK
V_n s_{n-1}\dots s_2s_1C_J.\end{equation} Let $\xi$ be a nonzero
element of $E_J$. Then using formula (41) we have
\begin{equation}\xi=\underset{u\in V_1}{\sum}a_{1,u}uC_J+ \underset{u\in V_2}{\sum}a_{2,u}us_1C_J+
\dots +\underset{u\in V_n}{\sum}a_{n,u}u s_{n-1}\dots
s_2s_1C_J,\end{equation} where $a_{i,u}\in \mathfrak K$ and only
finitely many of them are nonzero.

By Lemma 2.3 we can assume each $u$ appears in $\xi$ is
$\Delta_{w_Js_1s_2\dots s_{i-1}}$-regular wherever $a_{i,u}$ is
nonzero. Otherwise we replace $\xi$ by suitable $y\xi$. In this case
we also say that $\xi$ is regular.  Also we may choose $m$
sufficiently large so that all $n_i=n_{s_i}$ are in $G_{q^m}$.

For convenience, we can just write $s_ixwC_J$ instead of $n_ixwC_J$
for $x\in U$. With the assumption of $\xi$, using subsection 2.4 the
following results (a)-(c) are not difficult to see.

\medskip
(a) When $i=j-1$, we get $s_jus_{j-1}\dots s_2s_1C_J\in Us_j\dots
s_2s_1C_J$.

\medskip

(b) When $i=j$, we get $s_jus_j\dots s_2s_1C_J\in Us_j\dots
s_2s_1C_J$.

\medskip

(c) When $i\ne j-1$ and $i\ne j$, then $s_jus_i\dots
s_2s_1C_J=(x-y)s_i\dots s_2s_1C_J$ for some $x,y\in U$.

\medskip
Set
$$M'=\fK G\xi,\quad \mathbb{X}=\underset{x\in U_{q^m}}{\sum}x,
\quad \phi_{i,h}=\underset{u\in V_i}{\sum}h_{i,u}\ \ \text{for}\
h=a,\  a',\  b.$$ Using (a), (b), (c) we have
$$\mathbb{X}s_n\xi =
\phi_{n,a}\mathbb{X} s_n\dots s_2s_1C_J\in M'$$ If $\phi_{n,a} \ne
0$, by Lemma 2.9 then $s_n\dots s_2s_1C_J \in M'$ and $M'=E_J$.

Using (a), (b), (c) we have

$$\mathbb{X}s_j\xi = (\phi_{j,a}+\phi_{j+1,a})\mathbb{X} s_j\dots s_2s_1C_J \in M'$$

If $\phi_{j,a}+\phi_{j+1,a}\ne 0$, also by Lemma 2.9, then $s_j\dots
s_2s_1C_J \in M$ which implies $M'=E_J$. Therefore we know that if
$\phi_{i,a}\ne 0$ for some integer $i$, then $M'= E_J$.

To deal with the case  $\phi_{i,a} = 0$ for all $i$. We need the
following two lemmas.

\medskip

{\bf Lemma 5.2 \ } Let $I$ be any subset of $S$ and $u\in U_{w_I}$.
If $u$ is not the neutral element $1$, then the sum of all
coefficients of $w_IuC_I$ in terms the basis $zC_I,\ z\in U$, is 0.

\medskip

Proof. The argument of [Y, Lemma 2.5] can also be used here and
the proof is similar. \qed

\medskip

{\bf Lemma 5.3 \ } Let $x\in U_{q^m}$ and fix the element
$w^*=s_i\dots s_2s_1$ for some integer $i=1,2,\dots, n-1$. We write
$$ w_Jxw^*C_J =\underset{u\in V_1}{\sum}b_{1,u}uC_J+ \underset{u\in V_2}{\sum}b_{2,u}us_1C_J+
\dots +\underset{u\in V_n}{\sum}b_{n,u}u s_{n-1}\dots s_2s_1C_J$$
then we get $\phi_{1,b}+\phi_{2,b}=0$.

\medskip

Proof. Firstly, we compute the form of $s_jus_k\dots s_2s_1C_J $,
where $j\ge 2, k\ge 1$. Using subsection 2.4, we have the following
results.

\medskip

(d) When $j< k$, then $s_jus_k\dots s_2s_1C_J \in \mathfrak
KUs_k\dots s_2s_1C_J $. (The specific form depends on whether
$u_{\alpha_j}$ is the neutral element $1$.)

\medskip

(e) When $j=k$, if $u_{\alpha_j}=1$, then $s_jus_k\dots s_2s_1C_J =
u's_{k-1}\dots s_2s_1C_J$ for some $u'\in U$. If $u_{\alpha_j}\ne
1$, then $s_j us_k\dots s_2s_1C_J = u''s_k\dots s_2s_1C_J$ for some
$u''\in U$.

\medskip

(f) When $j= k+1$, then $s_jus_k\dots s_2s_1C_J=x s_{k+1}s_k\dots
s_2s_1C_J$ for some $x \in U$.

\medskip

(g) When $j> k+1$, then $s_jus_k\dots s_2s_1C_J \in \mathfrak
KUs_k\dots s_2s_1C_J $.

\medskip

Using (d), (e), (f), (g), we can get that
$$s_2s_3\dots s_n x w^*C_J\in \underset{2\le i \le
n}{\sum}\fK Us_i\dots s_2s_1C_J.$$

Since
$$\eta_S=C_J-s_1C_J+s_2s_1C_J-\dots+(-1)^n s_n\dots s_2s_1C_J$$ and
the image of $\eta_S$ is zero in $E_J$, we can see
\begin{equation}s_2s_3\dots s_n xw^*C_J = \underset{u\in
U}{\sum}b'_uu(C_J-s_1C_J) +\underset{3\le i \le
n}{\sum}\underset{u\in V_i}{\sum}b'_{i,u}u s_{i-1}\dots
s_2s_1C_J.\end{equation}

Let $K=\{3,4,\dots,n\}$ be the subset of $S$. Then it is easy to see
$w_J=w_K s_2s_3\dots s_n$.  Also using (d), (e), (f), (g), we have
that  $w_KuC_J \in \fK UC_J$, $w_Kus_1C_J \in \fK Us_1C_J$ and
\begin{equation}w_Kus_j\dots s_2s_1C_J\in \underset{2\le i \le
n}{\sum}\fK Us_{i}\dots s_2s_1C_J \ \ \text{for}\ j\ge
2.\end{equation} Therefore it is not difficult to see $w_JxwC_J$
satisfies our requirement in the lemma. The lemma is proved.\qed

\medskip

Using the two lemmas, we can deal with the case $\phi_{i,a}= 0$ for
all $i$.

Suppose $a_{1,u_1}\ne 0$ for some $u_1\in V_1$ and we consider
$w_J{u_1}^{-1}\xi$. Let $$w_J{u_1}^{-1}\xi= \underset{u\in
V_1}{\sum}a'_{1,u}uC_J+ \underset{u\in V_2}{\sum}a'_{2,u}us_1C_J+
\dots +\underset{u\in V_n}{\sum}a'_{n,u}u s_{n-1}\dots s_2s_1C_J.$$
Then by Lemma 5.2  and Lemma 5.3, we can see
$\phi_{1,a'}+\phi_{2,a'}=(-1)^{l(w_J)}a_{1,u_1}$. In this case, we
have show $M'= E_J$.

Suppose all $a_{1,u}$ are zero but $a_{2,u_2}\ne 0$ for some $u_2$.
Using (d) and (e) we have
$$s_1{u_2}^{-1}\xi=\underset{u\in V_1}{\sum}a''_{1,u}uC_J+ \underset{u\in V_2}{\sum}a''_{2,u}us_1C_J+
\dots +\underset{u\in V_n}{\sum}a''_{n,u}u s_{n-1}\dots s_2s_1C_J$$
and $a''_{1,1}=a_{2,u_2}$ is nonzero which implies $M'= E_J$.

Now suppose there exists an integer $k$ such that all $a_{j,u}=0$ for
$j<k$ but $a_{k,u_k}\ne 0$ for some $u_k$. Also using (d) and (e) we
have
$$s_{k-1}{u_k}^{-1}\xi= \underset{u\in V_{k-1}}{\sum}a^*_{k-1,u}us_{k-1}\dots s_2s_1C_J+
\dots +\underset{u\in V_n}{\sum}a^*_{n,u}u s_{n-1}\dots s_2s_1C_J$$
Thus we can deal with this case by induction and prove that $M'=
E_J$.

Therefore $E_J$ is irreducible. The theorem is proved. By symmetry,
$E_I$ is irreducible for $I=\{1,2,\dots, n-1\}$ a proper subset of
$S$.

\bigskip

{\bf 5.4}\ \ Now we consider the general case. Fix an integer $i'$
such that $2\le i' \le n-1$ and set $s=s_{i'}$. Let $J=S\setminus
\{s\}$ be a proper subset of $S$. In this section we prove $E_J$ is
irreducible.

As the notation in 4.4, set $m= i'-1$ and $l=n- i'$ and denote by
$r_j=s_{i'-j}$ for $j=1,2,\dots,m$, $t_k=s_{i'+k}$ for
$k=1,2,\dots,l$. The Dynkin diagram of $W$ is as follows:
\bigskip

\centerline{\begin{tabular}{lll}
&\xymatrix{\circ_{r_m}\ar@{-}[r]&\circ\cdots\cdots\circ_{r_2}\ar@{-}[r]&\circ_{r_1}\ar@{-}[r]&\circ_{s}\ar@{-}[r]&\circ_{t_1}\ar@{-}[r]&\circ\cdots\cdots\circ_{t_{l-1}}\ar@{-}[r]&\circ_{t_l}}
\end{tabular}}

\bigskip

As the notation in 2.5, Let $Y_J=\{z\in X_J\mid R(zw_J)=J\}$. In our
special case that $J$ is a proper subset of $S$, we have
$Y_J=X_J\backslash \{h_J\}$ where $h_J=w_0w_J$ is the longest
element in $X_J$.

Let $Y= Y_{\{s\}}$ be the subset of the Weyl group $W$ which has a
description in Lemma 4.5.  Thus we have $Y_J=(Ys \backslash
\{h_J\})\cup \{e\}$ where $e$ is the neutral element in $W$. Denote
by $d_J$ the unique element whose length is maximal in $Y_J$.

We show that $E_J$ is generated by any nonzero element in $E_J$ so
that $E_J$ is irreducible. Let $\xi$ be a nonzero element of $E_J$.
We denote by $V_w= U_{w_Jw^{-1}}$. Using Lemma 2.7 we have
\begin{equation}\xi=\underset{w\in Y_J}{\sum} \underset{x\in
V_w}{\sum}a_{w,x}xwC_J, \quad a_{w,x}\in\fK \end{equation} and only
finitely many of the coefficients $a_{w,x}$ are nonzero.

By Lemma 2.3 we can assume that $x$ appears in $\xi$ is
$\Delta_{w_Jw^{-1}}$-regular wherever $a_{w,x}$ is nonzero.
Otherwise, we replace $\xi$ by suitable $y\xi$. We can say that
$\xi$ is regular for simply in this case. We can choose $m$
sufficiently large so that all $x\in U_{q^m}$ if $a_{w,x}\ne 0$. We
also use the definition of regularlization given in subsection 4.4.

In subsection 4.4, we classify the elements of $Y_{\{s\}}$ by
different endings. We can do the same thing to $Y_J$ and let
$A'_k=A_ks\bigcap Y_J,\  B'_k=B_ks\bigcap Y_J,\  C'_k=C_ks\bigcap
Y_J,\  D'_k=D_ks\bigcap Y_J$. For convenience we write $A_k,\ B_k,\
C_k,\ D_k$ instead of $A'_k,\ B'_k,\ C'_k,\ D'_k$. Then $\{e\},\
A_k,\  B_k,\  C_k,\  D_k$ for $k=0,1,2, \dots , h$ give a partition
of $Y_J$ where $h$ is the maximal integer such that $A_h,\ B_h,\
C_h,\ D_h$ are nonzero.

Set
$$M'=\fK G\xi,\quad \mathbb{X}=\underset{x\in U_{q^m}}{\sum}x,
\quad \phi_{w,a}=\underset{x\in V_w}{\sum}a_{w,x}.$$ We also use the
notation
$$\sigma_i=r_ir_{i-1}\dots r_2r_1,\ \ \tau_j=t_jt_{j-1}\dots t_2t_1,\
\ w_c=r_m\dots r_2r_1t_l\dots t_2t_1s$$ as before.

\medskip

{\bf(Step 1)\ \ } We consider the element $R((r_m\dots
r_2r_1)(t_l\dots t_2t_1)s\xi)$ and denote this element by $\xi'$.
Using 2,4 (c), 2.4 (d), 2.4 (e) and the definition of
regularlization $R(\xi)$ of $\xi$, we have the following results.

\medskip

(a) We compute $sxA_kC_J$(resp. $sx'B_kC_J$). By 2.4 (e), we know
$sxwC_J=(y-z)wC_J$ for some $y,z\in U$ where $w\in A_k$. (resp.
$sx'w'C_J=(y'-z')wC_J$ for some $y',z'\in U$ where $w'\in B'_k$. )

\medskip

(b) Using 2.1 (a) and 2.4 (c), we know $sxC_kC_J\subseteq
UD_{k+1}C_J \cup Uh_JC_J$ where $h_J$ is the unique longest element
in $X_J$.

\medskip

(c) Using 2.4 (d), we know $sxD_kC_J\subseteq UD_kC_J$ for $k\ge 2$.

\medskip

If we use $\mathbb{X}$ multiply on $\xi'$, we can see the
coefficients $a_{w,x}$ will be killed for $w\in A_k\cup B_k$ by (a).
Using (b) and (c), it suffice to consider the form of $UD_kC_J$.

Let $s\tau_{j_k}\sigma_{i_k}s\dots \tau_{j_1}\sigma_{i_1}s$ be one
element in $D_k$ for $k\ge 2$.(Note $D_1$ just consist of $s$.)
Using 2.1(a), 2.4(c), 2.4(e), we can see

\medskip
(d)
$r_{i_k}x\tau_{{j_k}-1}\sigma_{{i_k}-1}s\tau_{j_k}\sigma_{i_k}s\dots
\tau_{j_1}\sigma_{i_1}sC_J$

$=(y-z)\tau_{{j_k}-1}\sigma_{{i_k}-1}s\tau_{j_k}\sigma_{i_k}s\dots
\tau_{j_1}\sigma_{i_1}sC_J$

for some $x, y, z\in U$.

\medskip

(e)
$t_{j_k}x'\tau_{{j_k}-1}\sigma_{{i_k}-1}s\tau_{j_k}\sigma_{i_k}s\dots
\tau_{j_1}\sigma_{i_1}sC_J$

$=(y'-z')\tau_{{j_k}-1}\sigma_{{i_k}-1}s\tau_{j_k}\sigma_{i_k}s\dots
\tau_{j_1}\sigma_{i_1}sC_J$

for some $x', y', z'\in U$.

\medskip

If we use $\mathbb{X}$ multiply on $\xi'$, we can see the
coefficients $a_{w,x}$ will be killed for $w\in C_k$ ($k\ge 1$) and
$w\in D_k$ ($k\ge 2$) by (d), (e).

On the other hand, We have following result by 2.4(c), 2.4(d).

\medskip
(f) $w_cxC_J\in Uw_cC_J$ and $w_cxsC_J\in Uw_cC_J$ where $x\in U$ is
$\alpha$-regular.

\medskip
Then we have
$$\mathbb{X}\xi' = (\phi_{e,a}+ \phi_{s,a})\mathbb{X}w_cC_J.$$
By Lemma 2.9, if $\phi_{e,a}+ \phi_{s,a} \ne 0$, then $w_cC_J\in M'$
which implies $M'=E_J$.

\medskip

{\bf(Step 2)\ \ } Let $r\in \{s,r_1,r_2,\dots,r_m,t_1,t_2,\dots,
r_l\}$ be a simple reflection. Using 2.4 (c), 2.4 (d), 2.4 (e) we
can compute the form of $rxwC_J$, where $w\in Y_J$.
\medskip

(a) When $l(rw)<l(w)$. If $x_{\alpha_r}=1$ then $rxwC_J=x'rwC_J$ for
some $x'\in U$. If $x_{\alpha_r}\ne 1$, then $rxwC_J=x''wC_J$ for
some $x''\in U$.

\medskip

(b) When $l(rww_J)>l(ww_J)$, then $rxwC_J=x^*rwC_J$ for some $x^*\in
U$. Note that when $d_J$ is the unique element whose length is
maximal in $Y_J$, then $rd_J C_J$ is the alternating sum of
$\{wC_J\mid w\in Y_J\}$.

\medskip

(c) When $l(rw)>l(w)$ but $l(rww_J)<l(ww_J)$. If $x_{\alpha_r}\ne
1$, then $rxwC_J=(y-z)wC_J$ for some $y, z \in U$. If $x_{\alpha_r}=
1$, then $rxwC_J=u wC_J$ for some $u\in U$.

\medskip

{\bf(Step 3)\ \ } As before $d_J$ is the unique element whose length
is maximal in $Y_J$. we compute the form of ${d_J}^{-1}xwC_J$, where
$w\in Y_J$ and $x$ is a general element in $V_w$. We give a claim:
\medskip

The element $ {d_J}^{-1}xwC_J$ has a expression form
\begin{equation} {d_J} ^{-1}xwC_J= \underset{w\in Y_J}{\sum}
\underset{x\in V_w}{\sum}b_{w,x}xwC_J + \underset{x\in U}{\sum}b_x
xh_JC_J\end{equation} such that the coefficients $b_{w_c,u}$ are
zero where $w_c=r_m\dots r_2r_1t_l\dots t_2t_1s$.
\medskip

Suppose the claim is not true, then the unique case may happen is
$w=d_J$. Indeed, if $w_cC_J$ appears in the expression form of $
{d_J}^{-1}xwC_J$ then $sw_cC_J$ has to appear in the expression form
of $s {d_J}^{-1}xwC_J$ by the results in step 2.

By (b) in step 2, we know $\sigma_{m-1}\tau_{l-1}sw_cC_J$ has to
appear in the expression form of $\sigma_{m-1}\tau_{l-1}
s\displaystyle {d_J}^{-1}xwC_J$. Using (c) in step 2, we know
$\sigma_{m-1}\tau_{l-1}sw_cC_J$ has to appear in the expression form
of $\sigma_m\tau_l s\displaystyle {d_J} ^{-1}xwC_J$. Using the
results in step 2 repeatedly, we can see if $w_cC_J$ appears in the
expression form of $\displaystyle {d_J} ^{-1}xwC_J$ then $d_J C_J$
has to appear in $xwC_J$ which implies $w=d_J$.

\medskip

Next we consider the expression form of ${d_J} ^{-1}x d_J C_J$. We
write $\displaystyle {d_J} ^{-1}x d_J C_J$ as the formula (46). We
pay attention on the number of $s$ occur in a reduced  expression of
$w\in Y_J$. Noting the number of $s$ occur in $w\in Y_J$ which we
denote by $l'(w)$ is independent of the reduced expression of $w$ so
it is well-defined.

Using the results in subsection 2.4, we compute the form of $sxwC_J$
and have following results.

\medskip

(a) If $l(sww_J)>l(ww_J)$ then $sxwC_J\in UswC_J$ and
$l'(sw)>l'(w)$.

\medskip

(b) If $l(sw)<l(w)$ and $x_{\alpha_s}\ne 1$ then $sxwC_J\in UwC_J$.
If $l(sw)>l(w)$ and $l(sww_J)<l(ww_J)$ then $sxwC_J\in \fK UwC_J$.

\medskip

(c) If $l(sw)<l(w)$ and $x_{\alpha_s}= 1$ then $sxwC_J\in UswC_J$
and $l'(sw)< l'(w)$.

\medskip

Suppose $w_cC_J$ appear in the expression form of ${d_J} ^{-1}x d_J
C_J$. Then by (a), (b), (c), the number of $s$ will be killed is
$l'(d_J)-1$ and only one $s$ leaves. Since $s$ being killed only in
the case of (c), the only $s$ leaves is the one in the right ending
of $d_J$. However we can not get $w_c$ in this case. We get a
contradiction. The claim is proved.

\medskip

{\bf(Step 4)\ \ } In this step, we give a lemma which is useful
later.

\medskip

{\bf Lemma 5.5 \ } Let $u\in U_{q^m}$ and $w$ be some element in
$Y_J$ such that $w\ne e$ and $w\ne w_c$. We write
\begin{equation} w_JuwC_J = \underset{w\in Y_J}{\sum} \underset{x\in
V_w}{\sum}b_{w,x}xwC_J\end{equation} and  then
$\phi_{e,b}+\phi_{s,b}=0.$

\medskip

Proof. Let $K= \{r_2, r_3,\dots, r_m, t_2,t_3,\dots, t_l\}$ be the
subset of $S$, then $w_J= w_K r_1r_2\dots r_m t_1t_2\dots t_l$. We
denote by $Z$ the subset of $Y_J$ consist of $A_k, B_k, C_k, D_k$
for $k\ge 2$. Using the result in step 2, we know that $s$ can not
be killed by $w_1=r_1r_2\dots r_m t_1t_2\dots t_l$. Then we have
the following results.

\medskip
(a) For $w\in Z$, we have $w_1xwC_J\in \underset{w\in
Z}{\sum}\mathfrak K UwC_J+ \mathfrak K Uh_JC_J$ where $h_J$ is the
longest element in $X_J$.

\medskip

(b) For $w\in Y_J\backslash Z$ but $w\ne e$ and $w\ne w_c$, we have
$w_1 xwC_J\in \underset{l(w)\ge 2}{\sum} \mathfrak K UwC_J.$

\medskip

(c) For $w=e$, we have $w_1 xC_J\in \mathfrak K UC_J.$

\medskip

(d) For $w=w_c$, we have $w_1 xwC_J\in \underset{l(w)\ge 2}{\sum}
\mathfrak K UwC_J \ \ \text{or}\ \ \mathfrak K UsC_J.$

\medskip

Using (a), (b), (c), (d) and noting that $h_JC_J$ is the alternating
sum of $\{wC_J\mid w\in Y_J\}$, we can have the following result.
Let $\zeta=w_1uwC_J$, if $w\ne e$ and $w\ne w_c$, then
\begin{equation}\zeta= \underset{x\in U}{\sum}b'_x x(C_J-sC_J)
+\underset{l(w)\ge 2}{\sum} \underset{x\in
V_w}{\sum}b'_{w,x}xwC_J.\end{equation}

Using the result in subsection 2.4, it is easy to see for any $r\in
K$, we get following results.

\medskip

(e) If $x_{\alpha_r}=1$, then $rx(C_J-sC_J)=x'(C_J-sC_J)$ for some
$x'\in U$.

\medskip

(f) If $x_{\alpha_r}\ne 1$, then $rx(C_J-sC_J)=(y-z)(C_J-sC_J)$ for
some $y,z\in U$.

\medskip

Using the result in step 2, then we have

\medskip

(g) If $l(w)\ge 2$, then $rxwC_J\in \underset{l(w)\ge 2}{\sum}
\mathfrak K UwC_J$.

\medskip

Using (e), (f), (g), $r\zeta$ also has the same expression of
formula (48) for $r\in K$. Therefore by $w_J= w_K r_1r_2\dots r_m
t_1t_2\dots t_l$, when  we write $ w_JuwC_J$ as the formula (47) we
have $\phi_{e,b}+\phi_{s,b}=0$. The lemma is proved.\qed

\medskip
{\bf(Step 5)\ \ }Using the results in step 3 and step 4, we compute
the form of $w_J{d_J}^{-1}xwC_J$, where $w\in Y_J$ and $x$ is a
general element in $V_w$.

We write \begin{equation}{d_J}^{-1}xwC_J= \underset{w\in Y_J}{\sum}
\underset{x\in V_w}{\sum}d_{w,x}xwC_J + \underset{x\in U}{\sum}d_x
xh_JC_J.\end{equation}

Using the results in step 2, we can see when $w\ne d_J$, then there
exists a expression form of ${d_J}^{-1}xwC_J$ like the formula (49)
such that all $d_{e,x}$ are zero. Using the results in step 3, we
can also have all $d_{w_c,x}$ are all zero, where $w_c=r_m\dots
r_2r_1t_l\dots t_2t_1s$.

We write \begin{equation}w_J{d_J}^{-1}xwC_J=\underset{w\in
Y_J}{\sum} \underset{x\in V_w}{\sum}h_{w,x}xwC_J.\end{equation} By
the Lemma 5.5 of step 4, we can see if $w\ne d_J$ then $\phi_{e,h}+
\phi_{s,h}=0$.

Next we deal with the case $w=d_J$. Using the result in step 2 and
step 3. We have the following results.

\medskip
(a) If $w= d_J$ and $x=1$, then $w_J{d_J} ^{-1}x d_J C_J=
(-1)^{l(w_J)}C_J$.

\medskip
(b) If $w=d_J$ and $x\ne 1$, we also write ${d_J}^{-1}x d_J C_J$ of
the formula (49). Then $d_{e,1}=0$ and $d_{w_1,x}$ are all zero. We
write $w_J{d_J} ^{-1}x d_J C_J$ as the formula of (50). By the Lemma
5.2 and Lemma 5.5, we have $\phi_{e,h}+ \phi_{s,h}=0$.

\medskip

In conclusion, If $w=d_J$ and $x=1$, then $w_J{d_J}^{-1}x d_J C_J=
(-1)^{l(w_J)}C_J$. For other case, $w_J{d_J} ^{-1}xwC_J$ has a
expression form as the formula (50) such that $\phi_{e,h}+
\phi_{s,h}=0$.

\medskip
{\bf(Step 6)\ \ } Now assume the coefficient $a_{d_J, x_0}$ of $\xi$
is nonzero for some $x_0\in U$. We write
$$w_J{d_J}^{-1}{x_0}^{-1}\xi= \underset{w\in Y_J}{\sum} \underset{x\in
V_w}{\sum}a'_{w,x}xw C_J.$$ Using the results in step 5, we can see
$\phi_{e,a'}+\phi_{s,a'}= (-1)^{l(w_J)}a_{d_J, x_0}C_J$.

By the results in step 1, in this case, we have $\mathfrak KG\xi= E_J$.

Now we assume all coefficient $a_{d_J, x}$ of $\xi$ are zero. Since
$\xi$ is nonzero, using the results in step 2 we can construct a new
element $\eta$ such that $\eta \in \mathfrak KG\xi$ and moreover,
$$\displaystyle \eta = \underset{w\in Y_J}{\sum} \underset{x\in
V_w}{\sum}a^*_{w,x}xw C_J$$ satisfies not all $a^*_{d_J, x}$ are
zero. Thus $\mathfrak KG\eta= E_J$ which implies $\mathfrak KG\xi=
E_J$

Therefore, $E_J$ is irreducible. Theorem 4.1 is proved.

\section{Induced modules of parabolic groups}

{\bf 6.1\ \ } In this section, we consider the induced modules of general parabolic subgroups. For any subset $I$ of $S$, $W_I$ is the standard parabolic subgroup of $W$. We let $\Delta_I:=\{\alpha\in \Delta \mid s_{\alpha}\in I\}$ and $\Phi_I:=\Phi\cap \sum_{\alpha\in \Delta_I}\mathbb{Z}\alpha$. Set $N_I$ the subgroup of $N$ containing $T$ such that $N_I/T=W_I$ and we denote by $P_I=BN_IB$ the standard parabolic subgroup of $G$ containing $B$. Then $P_I=\langle T,U_{\alpha}\mid \alpha \in \Phi^+\cup \Phi_I \rangle$.
We have Levi decomposition $P_I=U_I\rtimes L_I$, where
\begin{equation*} U_I:=\prod_{\alpha\in \Phi^+\backslash \Phi_I}U_\alpha\qquad L_I:=\langle T,U_{\alpha}\mid \alpha\in \Phi_I \rangle.\end{equation*}

\medskip

Let $\mathfrak K$ be a field. As the definition of spherical principal series representation in 2.4, we define the induced modules of parabolic groups. For a
one dimensional representation $\theta$ of $T$ over $\mathfrak K$,
let $\mathfrak K_\theta$ be the corresponding $\mathfrak KT$-module,
which will be regarded as $\mathfrak KP_I$-module through the natural
homomorphism $P_I\to T$.  We define the $\mathfrak KG$ module
$$M_I(\theta)=\mathfrak K G\otimes_{\mathfrak K P_I}\mathfrak K_\theta,$$
When $\theta$ is the trivial representation of $T$ over
$\mathfrak K$, we write $M_I(tr)$ for $M_I(\theta)$. When $I$ is empty, $P_I$ is the Borel subgroup $B$, then $M_I(tr)$ is just the spherical principal series representation of $G$.

\medskip

{\bf Lemma 6.2\ \ } Denote by $R(w)=\{s\in S\mid ws< w\}$ and $Z_I=\{w\in W \mid R(w)\subseteq S\backslash I \}$, then the $\fk G$-module $M_I(tr)$ is is the sum of all $\fk U_{w^{-1}}\,w1_{tr}$, where $w\in Z_I$. i.e. we have
$$M_I(tr)=\sum_{w\in Z_I}\mathfrak KU_{w^{-1}}\,w1_{tr}.$$

Proof.\ \ Using the Bruhat decomposition of $G$ we get $$M_I(tr)=\sum_{w\in W}\mathfrak KUw1_{tr}.$$

For $I\subseteq S$ and $w\in W$, there exists a unique decomposition $w=v\sigma$ such that $l(w)=l(v)+l(\sigma)$, $\sigma\in W_I$ and $R(v)\subseteq S\backslash I$.
By the definition of $M_I(tr)$, we can see $w1_{tr}=v1_{tr}$ which implies $$M_I(tr)=\sum_{w\in Z_I}\mathfrak KUw1_{tr}.$$

Numbering all positive roots in any order $\gamma_1,\gamma_2, \dots ,\gamma_r$, then $$U=U_{\gamma_1}U_{\gamma_2}\cdots U_{\gamma_r}.$$
Given a positive root $\gamma$, consider for which $\gamma$ we have $w^{-1}U_{\gamma}w\subseteq P_I$.
This implies $w^{-1}(\gamma)\in \Phi^+\cup \Phi_I$ by $w^{-1}U_{\gamma}w=U_{w^{-1}(\gamma)}$ and $P_I=\langle T,U_{\alpha}\mid \alpha \in \Phi^+\cup \Phi_I \rangle$. If $w^{-1}(\gamma)\in \Phi^+$, then we can see $U_{\gamma}\subseteq U'_{w^{-1}}$. Otherwise, $w^{-1}(\gamma)\in {\Phi_I}^-$ which implies $\gamma\in w({\Phi_I}^-)$. However for any simple root $\alpha\in \Delta_I$, since $R(w)\subseteq S\backslash I$, we can get $l(ws_{\alpha})=l(w)+l(s_{\alpha})$, where $s_{\alpha}$ is the simple reflection corresponded to $\alpha$. Then $w(\alpha)\in \Phi^+$ which implies $w({\Phi_I}^+)\subseteq \Phi^+$ and it contradicts that $\gamma$ is a positive root and $\gamma\in w({\Phi_I}^-)$.

Hence we can see only when $U_{\gamma}\subseteq U'_{w^{-1}}$, we have $w^{-1}U_{\gamma}w\subseteq P_I$. Therefore by 2.1 (c) we get
$$M_I(tr)=\sum_{w\in Z_I}\mathfrak KU_{w^{-1}}\,w1_{tr}.$$
The lemma is proved.\qed

\medskip
It is clearly that $M_I(tr)$ is a quotient module of $M(tr)$. The following theorem gives a explicit link between these induced modules and the spherical principal series representation.

\medskip
{\bf Theorem 6.3\ \ }For any subset $I\subseteq S$, we denote by $\displaystyle N_I(tr)=\sum_{s\in I} M(tr)_{\{s\}}$. Then we have $M_I(tr)=M(tr)/ N_I(tr)$.

Proof. \ \ Firstly we have a natural surjective homomorphism $\Theta: M(tr) \longrightarrow M_I(tr)$ such that $\Theta(g 1_{tr})=g 1_{tr}.$ It is obvious that $N_I(tr)\subseteq \text{Ker}\ \Theta$ by the definition $M(tr)_{\{s\}}=kUW(1-s)1_{tr}$ and $\displaystyle N_I(tr)=\sum_{s\in I} M(tr)_{\{s\}}$. We will show that
$\text{Ker}\ \Theta \subseteq N_I(tr)$ and then we get $M_I(tr)=M(tr)/ N_I(tr)$.

Let $\xi \in M(tr)$ be a element in $\text{Ker}\ \Theta$. We can write $\xi$ as
\begin{equation}\xi=\underset{w\in W}{\sum}\  \underset{x\in
U}{\sum}\ a_{w,x}xw1_{tr},\quad a_{w,x}\in\fk.\end{equation}

For $I\subseteq S$ and $w\in W$, there exists a unique decomposition $w=v\sigma$ such that $l(w)=l(v)+l(\sigma)$, $\sigma \in W_I$ and $R(v)\subseteq S\backslash I$. Use these notions, for any $v\in W$ such that $R(v)\subseteq S\backslash I$, we set
\begin{equation}I(v)=\{w\in W \mid w=v\sigma\  \text{such that} \ l(w)=l(v)+l(\sigma)\ \text{and}\ \sigma \in W_I \}.\end{equation}

For a fixed $v\in W$, the element $x\in U$ can be uniquely written as $x=yz$ such that $y\in U_{v^{-1}}$ and $z\in U'_{v^{-1}}$. For $y\in U_{v^{-1}}$, denote by
\begin{equation}v(y)=\{ x\in U\mid \ \text{for some}\  z\in U'_{v^{-1}}, \ x=yz\}.\end{equation}

For $\xi$ in the formula (51), $v\in W$ such that $R(v)\subseteq S\backslash I $ and $y\in U_{v^{-1}}$. By the notation (52) and (53) we let
\begin{equation}H(\xi,v, y)=\underset{w\in I(v)}{\sum}\  \underset{x\in
v(y)}{\sum}\ a_{w,x}xw1_{tr},\quad a_{w,x}\in\fk.\end{equation}
Thus it is easy to see that $H(\xi,v, y)$ is also in $\text{Ker}\ \Theta$.
By Lemma 6.2, we can get the sum of coefficients $a_{w,x}$ in $H(\xi,v, y)$ is zero, i.e, $\displaystyle \underset{w\in I(v)}{\sum}\  \underset{x\in v(y)}{\sum}\ a_{w,x}=0$.

Therefore it suffices to prove that $xw1_{tr}-yv1_{tr}$ is in $N_I(tr)$, where $w\in I(v)$ and $x\in
v(y)$. Then we can get that
$$H(\xi,v, y)=\underset{w\in I(v)}{\sum}\  \underset{x\in v(y)}{\sum}\ a_{w,x}(xw1_{tr}-yv1_{tr})$$
is also in $N_I(tr)$.

In $\fk G$-module $M(tr)$, we have $xv1_{tr}=yv1_{tr}$, where $x\in v(y)$. Since $w\in I(v)$, there exists a element $\sigma \in W_I$ such that $w=v\sigma$ and $l(w)=l(v)+l(\sigma)$. Therefore it is easy to see
$$xw1_{tr}-yv1_{tr}=x(w1_{tr}-v1_{tr})$$
is in the $\fk G$-module $N_I(tr)$. The theorem is proved.\qed

\medskip

Let $G$ be a connected reductive group over $\mathbb{F}_q$ such that its derived subgroup is of type $A_n$.
Suppose $I$ is a maximal proper subset of $S$ and then $P=P_I$ is a maximal parabolic subgroup of $G$. Using Theorem 4.1 and  Theorem 6.3 we can get that $M_I(tr)$ has a unique quotient module which is trivial and
a unique irreducible submodule.

\section {Steinberg Module}

N,Xi studied the infinite dimensional representations of $G$ by taking the direct limit of the finite dimensional representations of $G_{q^a}$, where $G_{q^a}$ is the $\mathbb{F}_{q^a}$-points of $G$.
Let $\mathscr{F}$ be the abelian category of all $\fk G$-modules. Then the Grothendieck group  $[\mathscr{F}]$ of $\mathscr{F}$ is generated by the set $\{[M]\mid M\in \mathscr{F}\}$ of isomorphic classes of $\fk G$-modules and the following relations: For short exact sequence $$0\longrightarrow M \longrightarrow N \longrightarrow L \longrightarrow 0$$ of $\fk G$-modules add the relation $[M]-[N]+[L]=0$.

\medskip

Let $A=\fk G$ and $A_a=\fk G_{q^a}$ for any positive integer $a$.  As in [X, 1.8], we can consider a category $\mathscr{F}_0$ of $A$-modules whose objects are those $A$-modules $M$ with a finite dimensional $A_i$-submodule $M_i$ for each $i$ such that $M$ is the union of all $M_i$ and for any positive integers $i,j$,  $M_i$ and $M_j$ are contained in $M_r$ whenever both $A_i$ and $A_j$ are
contained in $A_r$. Let $(M,M_i)$ and $(N,N_i)$ be two objects in  $\mathscr{F}_0$. The morphisms from $(M,M_i)$ to $(N,N_i)$ are just those homomorphisms of $A$-module from $M$ to $N$ such that $f(M_i)\subset N_i$ for all $i$. Clearly $\mathscr{F}_0$ is an abelian category.

\medskip

Let $P_I$ be a parabolic subgroup containing $B$. We denote by $P_{I,a}$ the $\mathbb{F}_{q^a}$-points of $P_I$. For two positive integers $a$ and $b$ such that $a \mid b$, we have a natural injective $G_{q^a}$-module homomorphism $\text{Ind}_{P_{I,a}}^{G_{q^a}} 1_{tr}\to \text{Ind}_{P_{I,b}}^{G_{q^b}} 1_{tr}$. This family of injections forms a direct system and the direct limit is $M_I(tr)$. Thus the object $(M_I(tr), \ \text{Ind}_{P_{I,a}}^{G_{q^a}} 1_{tr})$ is in the ableian category $\mathscr{F}_0$.

\medskip

{\bf Theorem 7.1\ \ }
Denote by $St$ the infinite dimensional Steinberg module, then we have \begin{equation}[St]=\sum_{I\subseteq S}(-1)^{|I|}[M_I(tr)]\end{equation} in the Grothendieck group $[\mathscr{F}]$.

Proof.\ \ For any positive integer $a$, denote by $St_a$ the Steinberg module of $\fk G_{q^a}$ which is
firstly introduced in [S]. Therefore we have
$$[St_a]=\sum_{I\subseteq S}(-1)^{|I|}[\text{Ind}_{P_{I,a}}^{G_{q^a}} 1_{tr}]$$ in the representation ring of $G_{q^a}$. Since the direct limit is an exact functor, we get the formula (55). The theorem is proved.\qed

\medskip

\section{Some questions }
Theorem 3.1 and Theorem 4.1 suggest the following conjecture.

\medskip

{\bf Conjecture 8.1\ \ } Let $G$ be a connected reductive group over
$\bar{\mathbb F}_q$ and $S$ the set of simple reflections of the
Weyl group of $G$. Assume that char\,$\fk
\ne$\,char\,$\bar{\mathbb{F}_q}$. Then the $\fk G$-module $E_J$
defined in [X, 2.6] (see also subsection 2.5 this paper) is
irreducible for any subset $J$ of $S$. Thus, in this case $M(tr)$
has $2^{|S|}$ composition factors, which are naturally ono-to-one
corresponding to the subsets of $S$.

\medskip

There is another evidence of this conjecture to be believable. Let $G$ be a connected reductive group over $\bar{\mathbb F}_q$ and $\fk$ be a field such that char\,$\fk \ne$\,char\,$\bar{\mathbb{F}_q}$. Now we assume this conjecture is ture, then using Theorem 6.3 we can deduce that the alternating sum expression (55) of infinite dimensional Steinberg module is true.

\medskip

When $\fk=\mathbb C$, Xi proved that $M(\theta)$ has only finitely many composition factors (see [X, Theorem 3.4]. Based on Xi's result and Theroems 3.1 and 4.1, it seems reasonable to suggest that $M(\theta)$ has only finitely many compostion factors for any character $\theta: T\to \fk^*$ provided that char\,$\fk\ne$\,char\,$\mathbb F_q$.

\medskip

It should also be interesting to study the composition factors of $M(\theta)$ provided that char\,$\fk=$\,char\,$\mathbb F_q$. In this case, X. Chen gave a necessary and sufficient condition of $M(\theta)$ to be irreducible (see [Ch]).

\bigskip

{\bf Acknowledgement.} I thank Professor Nanhua Xi for suggesting the topic of this paper and for  guidance.
The author also would like to thank Professor C.Bonnafe and Professor Ming Fang for helpful discussion and comments.

\end{document}